\documentclass[reqno]{amsart}

\usepackage{mathrsfs}
\usepackage{amsfonts}
\usepackage{amssymb}
\usepackage{xcolor}
\usepackage{caption}

\usepackage{algorithm}
\usepackage{algorithmic}

\usepackage[colorlinks,citecolor = red,
linkcolor=blue,hyperindex,CJKbookmarks]{hyperref}
\usepackage[hyphenbreaks]{breakurl}
\usepackage{cite}
\usepackage{graphicx}
\usepackage{multirow}

\usepackage{float}
\usepackage{subcaption}
\usepackage{cases}
\usepackage{listings}
\newtheorem{theorem}{Theorem}[section]
\newtheorem{lemma}[theorem]{Lemma}

\newtheorem{prop}[theorem]{Proposition}

\newtheorem{assumption}[theorem]{Assumption}

\theoremstyle{definition}
\newtheorem{definition}[theorem]{Definition}

\newtheorem{remark}[theorem]{Remark}

\newcommand{\ba}{\begin{array}}
	\newcommand{\ea}{\end{array}}
\newcommand{\bea}{\begin{eqnarray}}
	\newcommand{\eea}{\end{eqnarray}}

\newcommand{\R}{{\mathbb{R}}}

\numberwithin{equation}{section}

\allowdisplaybreaks

\begin{document}
	
\title[ADMM for nonseparable problems]{An Extended ADMM for 3-Block Nonconvex Nonseparable Problems with Applications}
	
\author{Zekun Liu}
\address{School of Mathematical Sciences, Shanghai Jiao Tong University, Shanghai 200240, China}
\email{sjtu\_lzk@sjtu.edu.cn}
	
%\author{Jinyan Fan*}
%\address{School of Mathematical Sciences, and MOE-LSC, Shanghai Jiao Tong University, Shanghai 200240,	China}
%\email{jyfan@sjtu.edu.cn}
%	
%\thanks{* Corresponding author}
%\thanks{The authors are supported by .}
	
	
	
\maketitle

\begin{abstract}
	We consider a 3-block Alternating Direction Method of Multipliers (ADMM) for solving nonconvex nonseparable problems with a linear constraint. Inspired by \cite[Sun, Toh and Yang, \textit{SIAM Journal on Optimization}, 25 (2015), pp.882-915]{wtwice}, the proposed ADMM follows the Block Coordinate Descent (BCD) cycle order $1\to 3\to 2\to 3$. We analyze its convergence based on the Kurdyka-Łojasiewicz property. We also discuss two useful extensions of the proposed ADMM with $2\to 3\to 1\to 3$ Gauss-Seidel BCD cycle order, and with adding a proximal term for more general nonseparable problems, respectively. Moreover, we make numerical experiments on two nonconvex problems: robust principal component analysis and nonnegative matrix completion. Results show the efficiency and outperformance of the proposed ADMM.
\end{abstract}

\section{Introduction}\label{sec:1}

In this paper, we consider the nonconvex nonseparable optimization problem:
\begin{equation}\label{equ:opt}
	\begin{aligned}
		\min_{X,Y,Z\in \R^{m\times n}} \;&\Phi(X)+\Psi(Y)+f(Z)+H(X,Y,Z)\\
		\mathrm{s.t.}\quad\quad &\mathcal{A}(X)+\mathcal{B}(Y)+\mathcal{C}(Z)=D,
	\end{aligned}
\end{equation}
where $\Phi: \R^{m\times n}\to\R\cup\{+\infty\}$ is proper lower semicontinuous, possibly nonsmooth and nonconvex, $\Psi: \R^{m\times n}\to\R\cup\{+\infty\}$ is convex and possibly nonsmooth, and $f: \R^{m\times n}\to\R$ is Lipschitz continuously differentiable with the modulus $L_f>0$ and possibly nonconvex, $H: \R^{m\times n}\times \R^{m\times n}\times \R^{m\times n} \to\R$ is smooth, $\mathcal{A}, \mathcal{B}, \mathcal{C}: \R^{m\times n}\to\R^{d_1\times d_2}$ are linear operators, and $D\in \R^{d_1\times d_2}$ is a given matrix. For ease of exposition, we assume $X$, $Y$, and $Z$ have the same dimensions here. A myriad of practical problems can be reduced to the special form of \eqref{equ:opt}, such as compressed sensing \cite{cs1,cs2,MMVADMM}, low-rank matrix sensing \cite{lowrankmatrix1,lowrankmatrix2,ADM} and robust principle component analysis \cite{RPCA1,RPCA2,RPCA3}.

Denote $\beta>0$ as the penalty parameter and $\Lambda \in \R^{d_1\times d_2}$ as the Lagrange multiplier. The augmented Lagrange function for \eqref{equ:opt} is defined as:
\begin{equation}\label{equ:alf}
	\begin{aligned}\mathcal{L}_\beta&(X,Y,Z,\Lambda):=\Phi(X)+\Psi(Y)+f(Z)+H(X,Y,Z)\\ &-\left\langle \Lambda,\mathcal{A}(X)+\mathcal{B}(Y)+\mathcal{C}(Z)-D\right\rangle+\frac{\beta}{2}\left\|\mathcal{A}(X)+\mathcal{B}(Y)+\mathcal{C}(Z)-D\right\|_F^2.\end{aligned}
\end{equation}
As an efficient first-order algorithm for separable problems derived from the Douglas–Rachford splitting method \cite{DRS1,DRS2}, the classic ADMM \cite{ClaADMM1,ClaADMM2} follows the usual Block Coordinate Descent (BCD) cycle order $1\to 2\to 3$ to solve \eqref{equ:opt}:
\begin{align}\label{equ:Claadmmscheme}
\left\{\begin{aligned}
	&X^{k+1}\in\mathop{\arg\min}\limits_{X} \mathcal{L}_{\beta}\left(X,Y^k,Z^k,\Lambda^k\right),  \\
	&Y^{k+1}\in\mathop{\arg\min}\limits_{Y} \mathcal{L}_{\beta}\left(X^{k+1},Y,Z^{k},\Lambda^k\right),  \\
	&Z^{k+1}=\mathop{\arg\min}\limits_{Z} \mathcal{L}_{\beta}\left(X^{k+1},Y^{k+1},Z,\Lambda^k\right),  \\
	&\Lambda^{k+1}=\Lambda^k-\beta\left(\mathcal{A}\left(X^{k+1}\right)+\mathcal{B}\left(Y^{k+1}\right)+\mathcal{C}\left(Z^{k+1}\right)-D\right). 
\end{aligned}\right.
\end{align}

When $f\equiv0$, $H\equiv0$ and $\Phi$, $\Psi$ are proper closed convex functions, \eqref{equ:Claadmmscheme} reduces to the classic 2-block ADMM which is designed for 2-block separable convex problems. A large body of literature has studied its convergence \cite{cvx2admm1,cvx2admm2,cvx2admm3,cvx2admm4}. Besides, considering $\Phi$ to be nonconvex, quantities of researches have established the convergence by the powerful Kurdyka-Łojasiewicz (KL) inequality \cite{routine3,ncvx2admm}. Furthermore, adding a nonseparable term $H$, Gao and Zhang \cite{cvxnsep2admm1} established the convergence of the proximal ADMM when $\Phi$ is convex, $\Psi$ is strongly convex and $\nabla H$ is Lipschitz continuous. Guo, Han, and Wu \cite{routine2} proved the global convergence of the classic ADMM under assumptions that $\Phi$ is semiconvex, $\nabla \Psi$ is Lipschitz continuous, and $\nabla H$ is Lipschitz continuous on bounded subsets based on the KL property.

Consider multi-block separable ADMM, i.e., $\Phi$, $\Psi$ and $f$ are all non-zero, and $H\equiv 0$. For this case, even though all objective functions are convex, the direct extension of classic ADMM to 3-block problems can also be divergent \cite{div3admm}. Han and Yuan \cite{strcvx3admm} first established the global convergence of the direct extension of classic ADMM to multi-block cases with all of the objective functions to be strongly convex and the penalty factor to be smaller than a threshold. Based on the KL property, Guo et al. \cite{ncvx3admm1} obtained the global convergence of the nonconvex classic ADMM for multi-block cases under assumptions that each $A_i$ is of full column rank and the penalty parameter is restricted to a certain range. Wang, Cao, and Xu \cite{BADMM} analyzed the convergence of the multi-block ADMM with the Bregman distance. Inspired by Sun, Toh and Yang \cite{wtwice}, Zhang et al. \cite{ztwice} proposed an extended proximal ADMM for nonconvex 3-block problems with the special BCD cycle order $1\to 3\to 2\to 3$, and established the global convergence with the help of the KL property. 

As for multi-block ADMM with coupled variables, it should be pointed out that the convergence of multi-block ADMM for \eqref{equ:opt} is still an open problem \cite{ncvxnsmonsepadmm}. Hong, Luo, and Razaviyayn \cite{ncvxnsmonsepadmm} considered the nonconvex sharing problem which is a special form of \eqref{equ:opt}, proving the convergence of the classic ADMM, together with some extensions using various block selection rules. Wang, Yin, and Zeng \cite{ncvxnsepadmm} considered the multi-block ADMM for a more general nonconvex nonsmooth case which can be nonseparable, and obtained the convergence with the objective function being continuous and each separable terms $f_i$ satisfying the prox-regularity. 

In this paper, we focus on the general nonconvex, nonsmooth and nonseparable problem \eqref{equ:opt} under mild assumptions on the coupled term $H$. Inspired by the works in \cite{wtwice,ztwice}, we also use the BCD cycle order $1\to 3\to 2\to 3$ to ensure the convergence of ADMM and call it ADMMn. The iteration scheme of ADMMn is given as follows:
\begin{subnumcases}
	{}
	X^{k+1}\in\mathop{\arg\min}\limits_{X} \mathcal{L}_{\beta}\left(X,Y^k,Z^k,\Lambda^k\right), \label{equ:admmscheme1} \\
	Z^{k+\frac12}=\mathop{\arg\min}\limits_{Z} \mathcal{L}_{\beta}\left(X^{k+1},Y^k,Z,\Lambda^k\right), \label{equ:admmscheme2} \\
	Y^{k+1}\in\mathop{\arg\min}\limits_{Y} \mathcal{L}_{\beta}\left(X^{k+1},Y,Z^{k+\frac12},\Lambda^k\right), \label{equ:admmscheme3} \\
	Z^{k+1}=\mathop{\arg\min}\limits_{Z} \mathcal{L}_{\beta}\left(X^{k+1},Y^{k+1},Z,\Lambda^k\right), \label{equ:admmscheme4} \\
	\Lambda^{k+1}=\Lambda^k-\beta\left(\mathcal{A}\left(X^{k+1}\right)+\mathcal{B}\left(Y^{k+1}\right)+\mathcal{C}\left(Z^{k+1}\right)-D\right). \label{equ:admmscheme5}
\end{subnumcases}
The main difference between ADMMn and the methods proposed in \cite{ztwice,wtwice} is that ADMMn is designed for the nonseparable problem \eqref{equ:opt} and guaranteed to be globally convergent, whereas others are not. To the best of our knowledge, there has been no previous research applying the extended ADMM with the BCD cycle order $1\to 3\to 2\to 3$ to the general 3-block nonconvex nonseparable problem \eqref{equ:opt}. Our main contributions are as follows: (1) Establish the convergence of ADMMn under mild assumptions. (2) Discuss two modifications of ADMMn which are both useful in applications, and prove their convergence. The first modification is swapping the update order of $X$ and $Y$ in ADMMn. The swapped BCD cycle order $2\to 3\to 1\to 3$ is helpful when $\Psi \equiv 0$ in \eqref{equ:opt}, since ADMMn can still update $Z$ twice instead of degenerating into the classic 2-block ADMM. Besides, in our pre-experiments, ADMMn and pADMMz \cite{ztwice} involved in Section \ref{sec:4} with BCD cycle order of 2 before 1 showed a better stability and performance than those with BCD cycle order of 1 before 2. We conjecture the reason behind this phenomenon is that these two ADMM methods with BCD cycle order of 2 before 1 are more likely to generate bounded sequences, while the conditions for those with BCD cycle order of 1 before 2 to generate bounded sequences are stricter and hence are difficult to be satisfied sometimes. It implies that swapping 1 and 2 in the BCD cycle order of ADMM methods may be helpful in some applications. The second modification is adding proximal terms in the update of $X$ or $Y$. With proximal terms, ADMMn can deal with more general problems which do not satisfy (a2) and (a6) in Assumption \ref{ass:converge}. (3) Give a brand new model of the nonnegative matrix completion problem in numerical experiments, which demonstrates how to employ ADMMn in practice.

The structure of this paper is as follows. In Section \ref{sec:2}, we give some notations and recall some definitions and well-known results used in the convergence analysis. The global convergence of ADMMn  is constructed in Section \ref{sec:3}. Some numerical experiments are presented to validate the performance of ADMMn in Section \ref{sec:4}. Section \ref{sec:5} concludes this paper.

\section{Preliminaries}\label{sec:2}

We list some notations and definitions which will be used in the analysis.

\textbf{Notations.} $\left \| \cdot \right \|_F$ represents the Frobenius-norm of a matrix. $\left \| \cdot \right \|$ is the spectral norm of a matrix or linear operator. $\left \langle \cdot,\cdot \right \rangle$ denotes the standard inner product of two matrices. $A\ge 0$ means all elements of $A$ are nonnegative. $A\succeq B$ implies $A-B$ is semi-positive definite. $A\odot B$ ($A\oslash B$) denotes the element-wise product (division) of two matrices of equal sizes. The adjoint operator of $\mathcal{A}$ is denoted as $\mathcal{A}^*$. The identity operator is $\mathcal{I}$ and the identity matrix is $I$. We write $\mathrm{dist}(X,A):=\inf_{Y\in A}\left\|X-Y \right\| $ to denote the distance between a point $X$ and a nonempty set $A$. For a function $f:\R^{m\times n}\to \R$ and two constants $a<b$, define $\left[a<f<b\right]:=\left\lbrace X\in \R^{m\times n}:a<f(X)<b \right\rbrace $.

Recall the subdifferential of a convex function $f$ \cite{book3} is defined as
\begin{equation}\label{def:cvxsubf}
	\partial f(X):=\left\{G\in \mathbb{R}^{m\times n}: f(Y)\ge f(X)+\langle  G,Y-X \rangle ,\; \forall Y\in \mathrm{dom}f \right\}.
\end{equation}
Moreover, the generalized Fermat's rule \cite{book4} implies that $0\in \partial f(X)$ is a necessary condition for $X$ to be a local minimizer of $f$. We call $X$ a critical point of $f$ if it satisfies $0\in \partial f(X)$, and the set of all critical points of $f$ is denoted as $\mathrm{crit}f$. 

\begin{definition}\label{def:critical}
	Call $W^{\star}:=(X^{\star},Y^{\star},Z^{\star},\Lambda^{\star})$ a critical point of the augmented Lagrange function $\mathcal{L}_{\beta}$ \eqref{equ:alf}, if it satisfies
	\begin{subnumcases}
		{}
		\mathcal{A}^*(\Lambda^{\star})-\nabla_X H(X^{\star},Y^{\star},Z^{\star})\in \partial\Phi(X^{\star}) \label{equ:1optadmm1} \\
		\mathcal{B}^*(\Lambda^{\star})-\nabla_Y H(X^{\star},Y^{\star},Z^{\star})\in \partial\Psi(Y^{\star}) \label{equ:1optadmm2} \\
		\mathcal{C}^*(\Lambda^{\star})-\nabla_Z H(X^{\star},Y^{\star},Z^{\star})= \nabla f(Z^{\star}) \label{equ:1optadmm3} \\
		\mathcal{A}(X^{\star})+\mathcal{B}(Y^{\star})+\mathcal{C}(Z^{\star})-D=0 \label{equ:1optadmm4}
	\end{subnumcases}
\end{definition}

Now we review the Kurdyka–Łojasiewicz (KL) property. 
Denote $\Phi_{\eta} (\eta>0)$ as the class of concave functions $\varphi: [0,\eta)\to \R_{+}$ satisfying: (1) $\varphi(0)=0$; (2) $\varphi$ is continuously differentiable on $(0,\eta)$ and continuous at 0; (3) $\varphi^{\prime}(x)>0$ for all $x\in (0,\eta)$.

\begin{definition}[KL property, \cite{routine1,KL1}]\label{def:KL}
	Let $f$ be a proper lower semicontinuous function.
	$f$ is said to have the KL property at $X^{\star} \in \mathrm{dom}\partial f:=\left\lbrace X\in \mathbb{R}^{m\times n}: \partial f(X) \ne \emptyset \right\rbrace $ if there exists an $\eta \in (0,+\infty]$, a neighborhood $U$ of $X^{\star}$, and a function $\varphi \in \Phi_{\eta}$ such that for all $X\in U\cap \left[f(X^{\star})<f<f(X^{\star})+\eta \right]$, the KL inequality holds:
	\begin{equation*}
		\varphi^{\prime}(f(X)-f(X^{\star}))\cdot \mathrm{dist}(0,\partial f(X))\geq1.
	\end{equation*}
	Call $f$ a KL function if it satisfies the KL property at each point of $\mathrm{dom}\partial f$.
\end{definition}

\begin{prop}[Uniformized KL property, \cite{limitset}]\label{pro:UKL}
	Let $f$ be a proper lower semicontinuous function and $\Omega$ be a compact set. Suppose $f\equiv f^{\star}$ is a constant on $\Omega$ and satisfies the KL property at each point of $\Omega$. Then, there exist $\eta,\sigma>0$ and $\varphi \in \Phi_{\eta}$ such that for all $X\in \left\{X\in\mathbb{R}^{m\times n}:\operatorname{dist}(X,\Omega)<\sigma \right\}\cap \left[f^{\star}<f<f^{\star}+\eta \right] $, it holds that
	\begin{equation*}
		\varphi^{\prime}(f(X)-f^{\star})\cdot \mathrm{dist}(0,\partial f(X))\geq1.
	\end{equation*}
\end{prop}

%\begin{pro}\label{pro:inequal}
%	We list some simple but useful inequalities:
%	\begin{itemize}
%		\item[(I1)] $\left\|\frac{\sum_{i=1}^n X_i}{n}\right\|_F^2\le \frac{1}{n}\sum_{i=1}^n \|X_i\|_F^2$.
%		\item[(I2)] $\sqrt{ab}\le \frac{a+b}{2}$ for any $a,b \in \mathbb{R}$.
%		\item[(I3)] $\left | \|A\|_F-\|B\|_F \right | \le \|A+B\|_F \le \|A\|_F+\|B\|_F$ for any $A,B \in \mathbb{R}^{m\times n}$.
%	\end{itemize}
%\end{pro}

%The following inequality for smooth functions is useful in the convergence analysis.
%\begin{lem}[\cite{book1}]\label{lem:GL-Lip}
%	Suppose $f:\mathbb{R}^{m\times n}\to \mathbb{R}$ is a continuously differentiable function with the modulus $L_f>0$, then, for any $X,Y\in \mathbb{R}^{m\times n}$, we have
%	\begin{displaymath}
%		|f(Y)-f(X)-\langle\nabla f(X),Y-X\rangle|\leq\frac {L_f}2\|Y-X\|_F^2.
%	\end{displaymath}
%\end{lem}

\section{Convergence analysis}\label{sec:3}

Before proving the convergence of ADMMn \eqref{equ:admmscheme1}-\eqref{equ:admmscheme5}, we first make the following assumptions.

\begin{assumption}\label{ass:converge}
 $\Psi,\mathcal{A},\mathcal{B},\mathcal{C},H$, and $\beta$ satisfy:
\begin{enumerate}
	\item[(a1)] $\Psi$ is continuous on its domain;
	\item[(a2)] $\mathcal{A}^*\mathcal{A}\succeq\mu_1\mathcal{I}$, $\mathcal{B}^*\mathcal{B}\succeq\mu_2\mathcal{I}$, $\mathcal{C}^*\mathcal{C}\succeq\mu_3\mathcal{I}$, and $\mathcal{C}\mathcal{C}^*\succeq\mu_4\mathcal{I}$ for some $\mu_1$, $\mu_2$, $\mu_3$, $\mu_4>0$;
	\item[(a3)] There exists $N(X)$ such that $\|\nabla_{X}H(X,Y_{1},Z_{1})-\nabla_{X}H(X,Y_{2},Z_{2})\|_F\leq N(X)(\|Y_{1}-Y_{2}\|_F+\|Z_{1}-Z_{2}\|_F)$ for any fixed $X$;
	\item[(a4)] $\nabla_Y H(X,Y,\cdot)$ is $L_1(X,Y)$-Lipschitz for any fixed $X$ and $Y$. $\nabla_Y H(X,\cdot,Z)$ is $L_2(X,Z)$-Lipschitz for any fixed $X$ and $Z$. $\nabla_Z H(X,Y,\cdot)$ is $L_3(X,Y)$-Lipschitz for any fixed $X$ and $Y$; 
	\item[(a5)] There exist $L_1,L_2,L_3>0,N>0$ such that $\sup_{k\in \mathbb{N}} L_1\left(X^{k+1},Y^{k+1}\right)\le L_1$, $\sup_{k\in \mathbb{N}} L_2\left(X^{k+1},Z^{k+\frac12}\right)\le L_2$, $\sup_{k\in \mathbb{N}} \left\lbrace L_3\left(X^{k},Y^{k}\right),L_3\left(X^{k+1},Y^{k}\right)\right\rbrace   \le L_3$, and $\sup_{k\in \mathbb{N}} N\left(X^{k}\right)\le N$;
	\item[(a6)] There exists $M>0$ such that $\|\nabla_{Z}H(X_{1},Y_{1},Z_{1})-\nabla_{Z}H(X_{2},Y_{2},Z_{2})\|_F\leq M(\|Y_{1}-Y_{2}\|_F+\|Z_{1}-Z_{2}\|_F)$;
	\item[(a7)] The penalty parameter $\beta$ in \eqref{equ:alf} satisfies
	\begin{align*}
		\hspace{3.5em}\beta&>\hat{\beta}:=\max \left\{ \frac{\mu_4 L_2+\sqrt{\mu_4^2 L_2^2+16\mu_2\mu_4(M+L_f)^2}}{2\mu_2\mu_4},\right. \\
		&\left. \frac{\mu_4 (L_f+L_3)+\sqrt{\mu_4^2 (L_f+L_3)^2+32\mu_3\mu_4(M+L_f)^2}}{2\mu_3\mu_4}, (M+L_f)\sqrt{\frac{\mu_4}{\mu_3}} \right\}.
	\end{align*}
\end{enumerate}
\end{assumption}

\begin{remark} Actually, Assumption \ref{ass:converge} (a6) implies $X$ and $Z$ are separable. One possible form of $H(X,Y,Z)$ can be $h_1(X,Y)+h_2(Y,Z)$. It is a mild assumption on the structure of the coupled term $H$ in practice. At the end of this section, we will also demonstrate the countermeasure when $H$ does not satisfy this assumption.
\end{remark}

We present the first-order optimality conditions for ADMMn here since it will be leveraged frequently:
\begin{subnumcases}
	{}
	0\in \partial\Phi(X^{k+1})+\nabla_X H(X^{k+1},Y^{k},Z^{k})-\mathcal{A}^*(\Lambda^{k})\notag \\
	\hspace{1.7em}+\beta \mathcal{A}^*\left(\mathcal{A}(X^{k+1})+\mathcal{B}(Y^{k})+\mathcal{C}(Z^{k})-D\right), \label{equ:2optadmm1} \\
	0=\nabla f(Z^{k+\frac{1}{2}})+\nabla_{Z}H(X^{k+1},Y^{k},Z^{k+\frac{1}{2}})-\mathcal{C}^*(\Lambda^{k})\notag \\
	\hspace{1.8em}+\beta \mathcal{C}^*\left(\mathcal{A}(X^{k+1})+\mathcal{B}(Y^{k})+\mathcal{C}(Z^{k+\frac12})-D\right), \label{equ:2optadmm2} \\
	0\in  \partial\Psi(Y^{k+1})+\nabla_Y H(X^{k+1},Y^{k+1},Z^{k+\frac12})-\mathcal{B}^*(\Lambda^{k}) \notag \\
	\hspace{1.7em}+\beta \mathcal{B}^*\left(\mathcal{A}(X^{k+1})+\mathcal{B}(Y^{k+1})+\mathcal{C}(Z^{k+\frac12})-D\right), \label{equ:2optadmm3} \\
	0=\nabla f(Z^{k+1})+\nabla_{Z}H(X^{k+1},Y^{k+1},Z^{k+1})-\mathcal{C}^*(\Lambda^{k})\notag \\
	\hspace{1.8em}+\beta \mathcal{C}^*\left(\mathcal{A}(X^{k+1})+\mathcal{B}(Y^{k+1})+\mathcal{C}(Z^{k+1})-D\right), \label{equ:2optadmm4} \\
	\Lambda^{k+1}=\Lambda^k-\beta\left(\mathcal{A}(X^{k+1})+\mathcal{B}(Y^{k+1})+\mathcal{C}(Z^{k+1})-D\right). \label{equ:2optadmm5}
\end{subnumcases}

For simplicity, denote $W^k:= \left(X^k,Y^k,Z^k,\Lambda^k\right)$ and $W^{\star}:= (X^{\star},Y^{\star},Z^{\star},\Lambda^{\star})$. We analyze the convergence by proving the following four lemmas as in \cite{routine1,ztwice,routine2,routine3}.

\begin{lemma}[Descent lemma]\label{lem:descent}
	Suppose $\left\{W^k\right\}_{k\in \mathbb{N}}$ is the sequence generated by ADMMn. If (a2), (a4), (a5), (a6) and (a7) in Assumption \ref{ass:converge} hold, then we have
	\begin{equation}\label{equ:descent}   \mathcal{L}_\beta(W^{k+1})\leq\mathcal{L}_\beta(W^k)-c_1\left\|Y^{k+1}-Y^k\right\|_F^2-c_2\left\|Z^{k+1}-Z^k\right\|_F^2,
	\end{equation}
	where $c_1:=\frac{\beta\mu_{2}-L_{2}}{2}-\frac{2(M+L_f)^{2}}{\beta\mu_{4}}>0$, $c_{2}:=\frac{\beta \mu_3 -L_f-L_3}{4}-\frac{2(M+L_f)^{2}}{\beta\mu_{4}}>0$.
\end{lemma}

\begin{proof}
	We analyze the following five descents one by one. First, 
	\begin{equation}\label{lem31-1}
		\begin{aligned}
			&\mathcal{L}_\beta(X^{k+1},Y^{k+1},Z^{k+1},\Lambda^{k+1})-\mathcal{L}_\beta(X^{k+1},Y^{k+1},Z^{k+1},\Lambda^k) \\
			=\hspace{0.6em} &-\left\langle\Lambda^{k+1}-\Lambda^k,\mathcal{A}(X^{k+1})+\mathcal{B}(Y^{k+1})+\mathcal{C}(Z^{k+1})-D\right\rangle   \\
			\overset{\eqref{equ:2optadmm5}}{=}& \frac1\beta\left\|\Lambda^{k+1}-\Lambda^k\right\|_F^2.
		\end{aligned}
	\end{equation}
	
	Second, we have
	\begin{equation}\label{lem31-2}
		\begin{aligned}
			&\mathcal{L}_\beta(X^{k+1},Y^{k+1},Z^{k+1},\Lambda^k)-\mathcal{L}_\beta(X^{k+1},Y^{k+1},Z^{k+\frac12},\Lambda^k) \\
			=\;& f(Z^{k+1})-f(Z^{k+\frac12})+H(X^{k+1},Y^{k+1},Z^{k+1})-H(X^{k+1},Y^{k+1},Z^{k+\frac12}) \\
			& -\left\langle \Lambda^k,\mathcal{C}(Z^{k+1}-Z^{k+\frac12})\right\rangle   +\frac{\beta}{2}\left\|\mathcal{A}(X^{k+1})+\mathcal{B}(Y^{k+1})+\mathcal{C}(Z^{k+1})-D\right \|_F^2 \\
			&-\frac{\beta}{2}\left\|\mathcal{A}(X^{k+1})+\mathcal{B}(Y^{k+1})+\mathcal{C}(Z^{k+\frac12})-D \right\|_F^2 .
		\end{aligned}
	\end{equation}
	Since $\nabla f(Z)$, $\nabla_Z H(X^{k+1},Y^{k+1},\cdot)$ is Lipschitz with $L_f$, $L_3(X^{k+1},Y^{k+1})$, respectively, it holds that 
	\begin{equation*}
		f(Z^{k+1})-f(Z^{k+\frac12})\le \left\langle\nabla f(Z^{k+1}),Z^{k+1}-Z^{k+\frac12}\right\rangle + \frac{L_f}{2}\left\|Z^{k+1}-Z^{k+\frac12} \right\|_F^2,
	\end{equation*}
	\begin{equation}\label{lem31-4}
		\begin{aligned}
			& H(X^{k+1},Y^{k+1},Z^{k+1})- H(X^{k+1},Y^{k+1},Z^{k+\frac12})\\
			\le & \left\langle\nabla_Z H(X^{k+1},Y^{k+1},Z^{k+1}),Z^{k+1}-Z^{k+\frac12}\right\rangle 
			+ \frac{L_3}{2}\left\|Z^{k+1}-Z^{k+\frac12} \right\|_F^2.
		\end{aligned}
	\end{equation}
	The polarization identity implies
	\begin{equation}\label{lem31-5}
		\begin{aligned}
			&\left\|\mathcal{A}(X^{k+1})\!+\!\mathcal{B}(Y^{k+1})\!+\!\mathcal{C}(Z^{k+1})\!-\!D\right \|_F^2\!-\!\left\|\mathcal{A}(X^{k+1})\!+\!\mathcal{B}(Y^{k+1})\!+\!\mathcal{C}(Z^{k+\frac12})\!-\!D\right \|_F^2\\
			= &\left\|\mathcal{C}(Z^{k+1}\!-\!Z^{k+\frac12}) \right\|_F^2\!+\! \left\langle \mathcal{A}(X^{k+1})\!+\!\mathcal{B}(Y^{k+1})\!+\!\mathcal{C}(Z^{k+\frac12})\!-\!D,\mathcal{C}(Z^{k+1}\!-\!Z^{k+\frac12})\right\rangle.
		\end{aligned}
	\end{equation}
	It follows from \eqref{equ:2optadmm4} that
	\begin{equation*}
		\begin{aligned}
			&\nabla f(Z^{k+1})+\nabla_Z H(X^{k+1},Y^{k+1},Z^{k+1})-\mathcal{C}^*(\Lambda^k) \\
			&+\beta \mathcal{C}^*\left(\mathcal{A}(X^{k+1})+\mathcal{B}(Y^{k+1})+\mathcal{C}(Z^{k+\frac12})-D\right)=-\beta \mathcal{C}^*\mathcal{C}\left(Z^{k+1}-Z^{k+\frac12}\right).
		\end{aligned}
	\end{equation*}
	Combining the above four expressions and \eqref{lem31-2}, we get
	\begin{equation}\label{lem31-7}
		\begin{aligned}
			&\mathcal{L}_\beta(X^{k+1},Y^{k+1},Z^{k+1},\Lambda^k)-\mathcal{L}_\beta(X^{k+1},Y^{k+1},Z^{k+\frac12},\Lambda^k) \\
			\le \hspace{0.4em} & \frac{L_f+L_3}{2}\left\|Z^{k+1}-Z^{k+\frac12} \right\|_F^2 - \frac{\beta}{2}\left\|\mathcal{C}\left( Z^{k+1}-Z^{k+\frac12}\right)  \right\|_F^2 \\
			\overset{\text{(a2)}}{\le}& \frac{L_f+L_3-\beta \mu_3}{2}\left\|Z^{k+1}-Z^{k+\frac12} \right\|_F^2.
		\end{aligned}
	\end{equation}
	
	Next, we have
	\begin{equation*}
		\begin{aligned}
			&\mathcal{L}_\beta(X^{k+1},Y^{k+1},Z^{k+\frac12},\Lambda^k)-\mathcal{L}_\beta(X^{k+1},Y^k,Z^{k+\frac12},\Lambda^k) \\
			=\;&\Psi(Y^{k+1})-\Psi(Y^k)+H(X^{k+1},Y^{k+1},Z^{k+\frac12})-H(X^{k+1},Y^k,Z^{k+\frac12})\\
			&-\left\langle \Lambda^k,\mathcal{B}(Y^{k+1}-Y^k)\right\rangle	+\frac{\beta}{2}\left\|\mathcal{A}(X^{k+1})+\mathcal{B}(Y^{k+1})+\mathcal{C}(Z^{k+\frac12})-D\right \|_F^2 \\
			&-\frac{\beta}{2}\left\|\mathcal{A}(X^{k+1})+\mathcal{B}(Y^k)+\mathcal{C}(Z^{k+\frac12})-D \right\|_F^2 .
		\end{aligned}
	\end{equation*}
	Because $\Psi$ is convex, by the definition of the subdifferential \eqref{def:cvxsubf}, we have
	\begin{equation}\label{lem31-9}
		\begin{aligned}
			\Psi(Y^{k+1})-\Psi(Y^k)\le \left\langle G^{k+1}, Y^{k+1}- Y^k \right\rangle, \; \forall  G^{k+1}\in \partial \Psi(Y^{k+1}).
		\end{aligned}
	\end{equation}
	Similar to \eqref{lem31-4} and \eqref{lem31-5}, it holds
	\begin{align*}
		& H(X^{k+1},Y^{k+1},Z^{k+\frac12})- H(X^{k+1},Y^{k},Z^{k+\frac12}) \\
		\le & \left\langle\nabla_Y H(X^{k+1},Y^{k+1},Z^{k+\frac12}),Y^{k+1}-Y^{k}\right\rangle+ \frac{L_2}{2}\left\|Y^{k+1}-Y^{k} \right\|_F^2.
	\end{align*}
	\begin{align*}
		&\left\|\mathcal{A}(X^{k+1})\!+\!\mathcal{B}(Y^{k+1})\!+\!\mathcal{C}(Z^{k+\frac12})\!-\! D\right \|_F^2\!-\!\left\|\mathcal{A}(X^{k+1})\!+\!\mathcal{B}(Y^{k})\!+\!\mathcal{C}(Z^{k+\frac12})\!-\! D\right \|_F^2 \\
		= &\left\|\mathcal{B}(Y^{k+1}\!-\! Y^{k})\right \|_F^2\!+\! \left\langle \mathcal{A}(X^{k+1})\!+\! \mathcal{B}(Y^{k})\!+\! \mathcal{C}(Z^{k+\frac12})\!-\! D,\mathcal{B}(Y^{k+1}\!-\! Y^{k})\right\rangle.
	\end{align*}
	From \eqref{equ:2optadmm3}, there exists $\Tilde{G}^{k+1}\in \partial \Psi(Y^{k+1})$ such that 
	\begin{align*}
		0=&\;\Tilde{G}^{k+1}+\nabla_Y H(X^{k+1},Y^{k+1},Z^{k+\frac12})-\mathcal{B}^*(\Lambda^{k})\\
		&+\beta \mathcal{B}^*\left(\mathcal{A}(X^{k+1})+\mathcal{B}(Y^{k+1})+\mathcal{C}(Z^{k+\frac12})-D\right).
	\end{align*}
	Let $G^{k+1}=\Tilde{G}^{k+1}$ in \eqref{lem31-9}, and combine the above five expressions to conclude that
	\begin{equation}\label{lem31-14}
		\begin{aligned}
			& \mathcal{L}_\beta(X^{k+1},Y^{k+1},Z^{k+\frac12},\Lambda^k)-\mathcal{L}_\beta(X^{k+1},Y^{k},Z^{k+\frac12},\Lambda^k) \\
			\le \hspace{0.4em}&\frac{L_2}{2}\left\|Y^{k+1}-Y^{k} \right\|_F^2 - \frac{\beta}{2}\left\|\mathcal{B}\left( Y^{k+1}-Y^{k}\right)  \right\|_F^2 \\
			\overset{\text{(a2)}}{\le}& \frac{L_2-\beta \mu_2}{2}\left\|Y^{k+1}-Y^{k} \right\|_F^2.
		\end{aligned}
	\end{equation}
	
	Analogous to \eqref{lem31-7}, we obtain
	\begin{equation}\label{lem31-15}
		\begin{aligned}
			&\mathcal{L}_\beta(X^{k+1},Y^{k},Z^{k+\frac12},\Lambda^k)-\mathcal{L}_\beta(X^{k+1},Y^{k},Z^{k},\Lambda^k) \\
			\le \; & \frac{L_f+L_3-\beta \mu_3}{2}\left\|Z^{k+\frac12}-Z^{k}\right \|_F^2.
		\end{aligned}
	\end{equation}
	
	At last, It follows from the update rule \eqref{equ:admmscheme1} that
	\begin{equation}\label{lem31-16}
		\begin{aligned}
			\mathcal{L}_\beta(X^{k+1},Y^{k},Z^{k},\Lambda^k)-\mathcal{L}_\beta(X^{k},Y^{k},Z^{k},\Lambda^k) \le 0.
		\end{aligned}
	\end{equation}
	
	Moreover, combining \eqref{equ:2optadmm4} and \eqref{equ:2optadmm5}, we have
	\begin{equation}\label{equ:2optadmm6}
		\mathcal{C}^*(\Lambda^{k+1})=\nabla f(Z^{k+1})+\nabla_Z H(X^{k+1},Y^{k+1},Z^{k+1}).
	\end{equation}
	Thus, 
	\begin{equation*}
		\mathcal{C}^*(\Lambda^{k+1}\!-\!\Lambda^{k})\!=\! \nabla f(Z^{k+1})\!-\!\nabla f(Z^{k})\!+\! \nabla_Z H(X^{k+1},Y^{k+1},Z^{k+1})\!-\! \nabla_Z H(X^{k},Y^{k},Z^{k}),
	\end{equation*}
	which implies that
	\begin{equation}\label{equ:lambda}
		\begin{aligned}
			\left\| \Lambda^{k+1}-\Lambda^{k} \right\|_F & \overset{\text{(a2)}}{\le} \frac{1}{\sqrt{\mu_4}} \left\|\mathcal{C}^*(\Lambda^{k+1}-\Lambda^{k})\right \|_F \\
			& \overset{\text{(a6)}}{\le} \frac{M+L_f}{\sqrt{\mu_4}}\left(\left\| Y^{k+1}-Y^{k}\right\|_F+\left\|Z^{k+1}-Z^{k} \right\|_F\right).
		\end{aligned}
	\end{equation}
	Hence,
	\begin{equation}\label{lem31-17}
		\left\| \Lambda^{k+1}-\Lambda^{k} \right\|_F^2 \le \frac{2(M+L_f)^2}{\mu_4}\left(\left\| Y^{k+1}-Y^{k}\|_F^2+\|Z^{k+1}-Z^{k} \right\|_F^2\right).
	\end{equation}
	
	Note from Assumption \ref{ass:converge} (a7) that $\frac{L_f+L_3-\beta \mu_3}{2}<0$, thus
	\begin{equation}\label{lem31-18}
		\begin{aligned}
			& \frac{L_f+L_3-\beta \mu_3}{2}\left(\left\|Z^{k+1}-Z^{k+\frac12}\right \|_F^2+\left\|Z^{k+\frac12}-Z^{k}\right \|_F^2\right) \\
			\le &\; \frac{L_f+L_3-\beta \mu_3}{4}\left\|Z^{k+1}-Z^{k} \right\|_F^2.
		\end{aligned}
	\end{equation}
	
	Combining \eqref{lem31-1}, \eqref{lem31-7}, \eqref{lem31-14}, \eqref{lem31-15}, \eqref{lem31-16}, \eqref{lem31-17}, and \eqref{lem31-18}, we obtain
	\begin{equation*}
		\begin{aligned}
			&\mathcal{L}_\beta(W^{k+1})-\mathcal{L}_\beta(W^k) \\
			\le \; & \left(\frac{L_2-\beta \mu_2}{2}+\frac{2(M+L_f)^2}{\beta \mu_4}\right) \left\|Y^{k+1}-Y^k\right\|_F^2 \\
			& + \left(\frac{L_f+L_3-\beta \mu_3}{4}+\frac{2(M+L_f)^2}{\beta \mu_4}\right) \left\|Z^{k+1}-Z^k\right\|_F^2 \\ =&
			-c_1\left\|Y^{k+1}-Y^k\right\|_F^2-c_2\left\|Z^{k+1}-Z^k\right\|_F^2,
		\end{aligned}
	\end{equation*}
	where $c_1:=\frac{\beta\mu_{2}-L_{2}}{2}-\frac{2(M+L_f)^{2}}{\beta\mu_{4}}>0$ and $c_{2}:=\frac{\beta \mu_3 -L_f-L_3}{4}-\frac{2(M+L_f)^{2}}{\beta\mu_{4}}>0$ by Assumption \ref{ass:converge} (a7). This completes the proof.
\end{proof}

\begin{lemma}[Square summability]\label{lem:square}
	Suppose the sequence $\left\{W^k\right\}_{k\in \mathbb{N}}$ generated by ADMMn is bounded. 
	Then, under the conditions of Lemma \ref{lem:descent} , we have
	\begin{equation*}
		\sum_{k=0}^{\infty}\left\|W^{k+1}-W^k\right\|_F^2<+\infty.
	\end{equation*}
\end{lemma}

\begin{proof}
	Since $\left\{W^k\right\}_{k\in \mathbb{N}}$ is bounded, it has a convergent subsequence $\left\{W^{k_j}\right\}_{j\in \mathbb{N}}$. Denote the limit point $W^{\star}:=\lim_{j\to \infty}W^{k_j}$. It is easy to check that $\mathcal{L}_{\beta}$ is lower semicontinuous, hence
	\begin{equation}\label{lem32-1}
		\mathcal{L}_\beta(W^{\star})\leq \liminf_{j\to\infty}\mathcal{L}_\beta(W^{k_j}).
	\end{equation}
	Thus $\left\{\mathcal{L}_\beta(W^{k_j})\right\}$ is bounded from below. From Lemma \ref{lem:descent} we know $\left\{\mathcal{L}_\beta(W^{k_j})\right\}$ is monotonically decreasing. So $\left\{\mathcal{L}_\beta(W^{k_j})\right\}$ is convergent. Combining $\left\{\mathcal{L}_\beta(W^{k})\right\}$ being monotonically decreasing, it holds that $\left\{\mathcal{L}_\beta(W^{k})\right\}$ is also convergent.
	
	By Lemma \ref{lem:descent}, we have
	\begin{equation*}
		c_1\left\|Y^{k+1}-Y^k\right\|_F^2+c_2\left\|Z^{k+1}-Z^k\right\|_F^2 \le
		\mathcal{L}_\beta(W^k)-\mathcal{L}_\beta(W^{k+1}).
	\end{equation*}
	Summing over $k$ from $0$ to $K$, then letting $K\to \infty$, combining $\left\{\mathcal{L}_\beta(W^{k})\right\}$ being convergent and \eqref{lem32-1}, we get
	\begin{equation*}
		c_1\sum_{k=0}^{\infty}\left\|Y^{k+1}-Y^k\right\|_F^2+c_2\sum_{k=0}^{\infty}\left\|Z^{k+1}-Z^k\right\|_F^2 \le \mathcal{L}_\beta(W^{0})-\mathcal{L}_\beta(W^{\star})<+\infty.
	\end{equation*}
	Therefore,
	\begin{equation*}
		\sum_{k=0}^{\infty}\left\|Y^{k+1}-Y^k\right\|_F^2<+\infty, \quad \sum_{k=0}^{\infty}\left\|Z^{k+1}-Z^k\right\|_F^2<+\infty.
	\end{equation*}
	Besides, \eqref{lem31-17} tells that 
	\begin{equation*}
		\sum_{k=0}^{\infty}\left\|\Lambda^{k+1}-\Lambda^k\right\|_F^2<+\infty.
	\end{equation*}
	Recall \eqref{equ:2optadmm5}
	\begin{equation}\nonumber
		\begin{aligned}
			\Lambda^{k+1}&=\Lambda^k-\beta\left(\mathcal{A}(X^{k+1})+\mathcal{B}(Y^{k+1})+\mathcal{C}(Z^{k+1})-D\right), \\
			\Lambda^{k}&=\Lambda^{k-1}-\beta\left(\mathcal{A}(X^{k})+\mathcal{B}(Y^{k})+\mathcal{C}(Z^{k})-D\right),
		\end{aligned}
	\end{equation}
	hence
	\begin{equation*}
		\Lambda^{k+1}-\Lambda^{k}=\Lambda^{k}-\Lambda^{k-1}-\beta\mathcal{A}(X^{k+1}-X^{k})-\beta\mathcal{B}(Y^{k+1}-Y^{k})-\beta\mathcal{C}(Z^{k+1}-Z^{k}).
	\end{equation*}
	Then it follows that
	\begin{equation*}
		\begin{aligned}
			&\beta \sqrt{\mu_1}\left\|X^{k+1}-X^{k}\right\|_F 
			\overset{\text{(a2)}}{\le} \left\|\beta\mathcal{A}(X^{k+1}-X^{k})\right\|_F \\
			= & \left\| (\Lambda^{k}-\Lambda^{k-1})-(\Lambda^{k+1}-\Lambda^{k})-\beta\mathcal{B}(Y^{k+1}-Y^{k})-\beta\mathcal{C}(Z^{k+1}-Z^{k})  \right\|_F \\
			\le & \left\|\Lambda^{k}-\Lambda^{k-1} \right \|_F+ \left\|\Lambda^{k+1}-\Lambda^{k} \right \|_F+ \left\|\beta\mathcal{B}(Y^{k+1}-Y^{k}) \right \|_F+ \left\|\beta\mathcal{C}(Z^{k+1}-Z^{k}) \right \|_F.
		\end{aligned}
	\end{equation*}
	Therefore, we get
	\begin{equation}\label{equ:X}
		\begin{aligned}
			\left\| X^{k+1}-X^{k}\right\| _F\le \frac{1}{\beta \sqrt{\mu_1}}&\left( \left\| \Lambda^{k}-\Lambda^{k-1}  \right\| _F+ \left\| \Lambda^{k+1}-\Lambda^{k}  \right\| _F \right.\\
			& \hspace{0.5em} \left. + \beta \|\mathcal{B}\|\left\| Y^{k+1}-Y^{k}\right\| _F+ \beta \|\mathcal{C}\| \left\| Z^{k+1}-Z^{k} \right\| _F\right),
		\end{aligned}
	\end{equation}
	which implies that $\sum_{k=0}^{\infty}\left\|X^{k+1}-X^{k}\right\|_F^2<+\infty$ by Jensen's inequality.
	Hence we obtain $\sum_{k=0}^{\infty}\left\|W^{k+1}-W^{k}\right\|_F^2<+\infty$.
\end{proof}

\begin{lemma}[Boundness of the subdifferential]\label{lem:boundsub}
	Suppose the sequence $\left\{W^k\right\}_{k\in \mathbb{N}}$ generated by ADMMn is bounded, and (a2), (a3), (a4), (a5), and (a6) in Assumption \ref{ass:converge} hold. Then there exists $\vartheta>0$ such that
	\begin{equation}\label{equ:boundsub}
		\mathrm{dist}\left( 0,\partial \mathcal{L}_{\beta}(W^{k+1})\right) \le \vartheta \left(\left\|Y^{k+1}-Y^{k}\right\|_F+\left\|Z^{k+1}-Z^{k}\right\|_F\right).
	\end{equation}
\end{lemma}

\begin{proof}
	By the definition of $\mathcal{L}_{\beta}$ in \eqref{equ:alf}, we have
	\begin{equation*}
		\begin{cases}
			\partial_X \mathcal{L}_{\beta}(W^{k+1})=\partial\Phi(X^{k+1})+\nabla_X H(X^{k+1},Y^{k+1},Z^{k+1})-\mathcal{A}^*(\Lambda^{k+1}) \\
			\hspace{7.2em} +\beta \mathcal{A}^*\left(\mathcal{A}(X^{k+1})+\mathcal{B}(Y^{k+1})+\mathcal{C}(Z^{k+1})-D\right),  \\
			\partial_Y \mathcal{L}_{\beta}(W^{k+1})=\partial\Psi(Y^{k+1})+\nabla_Y H(X^{k+1},Y^{k+1},Z^{k+1})-\mathcal{B}^*(\Lambda^{k+1})\\
		    \hspace{7.2em} +\beta \mathcal{B}^*\left(\mathcal{A}(X^{k+1})+\mathcal{B}(Y^{k+1})+\mathcal{C}(Z^{k+1})-D\right),  \\
			\nabla_Z \mathcal{L}_{\beta}(W^{k+1})=\nabla f(Z^{k+1})+\nabla_{Z}H(X^{k+1},Y^{k+1},Z^{k+1})-\mathcal{C}^*(\Lambda^{k+1})\\
			\hspace{7.4em} +\beta \mathcal{C}^*\left(\mathcal{A}(X^{k+1})+\mathcal{B}(Y^{k+1})+\mathcal{C}(Z^{k+1})-D\right),  \\
			\nabla_{\Lambda} \mathcal{L}_{\beta}(W^{k+1})=-\left(\mathcal{A}(X^{k+1})+\mathcal{B}(Y^{k+1})+\mathcal{C}(Z^{k+1})-D\right).  
		\end{cases}
	\end{equation*}
	Substituting \eqref{equ:2optadmm1}-\eqref{equ:2optadmm5} into it yields
	\begin{equation*}
		\begin{cases}
			\Xi_1^{k+1}\in \partial_X \mathcal{L}_{\beta}(W^{k+1}), \\
			\Xi_2^{k+1}\in \partial_Y \mathcal{L}_{\beta}(W^{k+1}), \\
			\Xi_3^{k+1}=\nabla_Z \mathcal{L}_{\beta}(W^{k+1}),\\
			\Xi_4^{k+1}=\nabla_{\Lambda} \mathcal{L}_{\beta}(W^{k+1}),
		\end{cases}
	\end{equation*}
	where
	\begin{equation*}
		\begin{cases}
			\Xi_1^{k+1}:=\mathcal{A}^*\left( \Lambda^{k}-\Lambda^{k+1}\right) +\beta \mathcal{A}^*\left [  \mathcal{B}\left( Y^{k+1}-Y^{k}\right) +\mathcal{C}\left( Z^{k+1}-Z^{k}\right) \right ] \\
			\hspace{4em} +\nabla_X H(X^{k+1},Y^{k+1},Z^{k+1})-\nabla_X H(X^{k+1},Y^{k},Z^{k}), \\
			\Xi_2^{k+1}:=\mathcal{B}^*\left( \Lambda^{k}-\Lambda^{k+1}\right) +\beta \mathcal{B}^*\left [  \mathcal{C}\left( Z^{k+1}-Z^{k+\frac12}\right) \right ] \\
			\hspace{4em} +\nabla_Y H(X^{k+1},Y^{k+1},Z^{k+1}) 
			-\nabla_Y H(X^{k+1},Y^{k+1},Z^{k+\frac12}), \\
			\Xi_3^{k+1}:=\mathcal{C}^*\left( \Lambda^{k}-\Lambda^{k+1}\right) ,\\
			\Xi_4^{k+1}:=\frac{1}{\beta}\left( \Lambda^{k+1}-\Lambda^{k}\right) .
		\end{cases}
	\end{equation*}
	We know $\left(\Xi_1^{k+1},\Xi_2^{k+1},\Xi_3^{k+1},\Xi_4^{k+1}\right)\in \partial \mathcal{L}_{\beta}(W^{k+1})$ (cf.  \cite[Proposition 2.1]{routine1}). Therefore, we conclude that
	\begin{equation}\label{lem33-4}
		\begin{aligned}
			&\; \mathrm{dist}\left( 0,\partial \mathcal{L}_{\beta}(W^{k+1})\right) 
			\le  \left\|\left(\Xi_1^{k+1},\Xi_2^{k+1},\Xi_3^{k+1},\Xi_4^{k+1}\right) \right\|_F \\
			\le& \left\|\Xi_1^{k+1}\right \|_F+\left\|\Xi_2^{k+1}\right \|_F+\left\|\Xi_3^{k+1} \right\|_F+\left\| \Xi_4^{k+1}\right\|_F\\
			\le &\left(\|\mathcal{A}^*\|+\|\mathcal{B}^*\|+\|\mathcal{C}^*\|+\frac{1}{\beta}\right)\left\| \Lambda^{k+1}-\Lambda^{k}\right \|_F+\left(\beta\|\mathcal{A}^*\|\|\mathcal{B}\| +N\right)\left\|Y^{k+1}-Y^{k} \right \|_F\\
			&+\left(\beta\|\mathcal{A}^*\|\|\mathcal{C}\|+N\right) \left\|Z^{k+1}-Z^{k}  \right\|_F +\left(\beta\|\mathcal{B}^*\|\|\mathcal{C}\|+L_1\right) \left\| Z^{k+1}-Z^{k+\frac12} \right\|_F.
		\end{aligned}
	\end{equation}
	Now we estimate $\left\|  Z^{k+1}-Z^{k+\frac12} \right\| _F$. By \eqref{equ:2optadmm2}, \eqref{equ:2optadmm4}, and \eqref{equ:2optadmm5}, we get
	\begin{equation*}
		\begin{aligned}
			0=&\nabla f(Z^{k+\frac{1}{2}})+\nabla_{Z}H(X^{k+1},Y^{k},Z^{k+\frac{1}{2}})-\mathcal{C}^*(\Lambda^{k+1}) \\
			& -\beta \mathcal{C}^*\left [\mathcal{B}\left( Y^{k+1}-Y^k\right) +\mathcal{C}\left( Z^{k+1}-Z^{k+\frac12}\right) \right ], \\
			0=&\nabla f(Z^{k+1})+\nabla_{Z}H(X^{k+1},Y^{k+1},Z^{k+1})-\mathcal{C}^*(\Lambda^{k+1}). 
		\end{aligned}
	\end{equation*}
	Then it follows
	\begin{equation}\label{lem33-5}
		\begin{aligned}
			&\left\|\nabla f(Z^{k+1})\!-\! \nabla f(Z^{k+\frac{1}{2}})\!+\! \nabla_{Z}H(X^{k+1},Y^{k+1},Z^{k+1})\!-\! \nabla_{Z}H(X^{k+1},Y^{k},Z^{k+\frac{1}{2}})\right\|_F \\
			= &\beta \left \| \mathcal{C}^*\left [\mathcal{B}\left( Y^{k+1}-Y^k\right) +\mathcal{C}\left( Z^{k+1}-Z^{k+\frac12}\right)  \right ] \right \|_F.
		\end{aligned}
	\end{equation}
	The LHS of \eqref{lem33-5} 
	\begin{equation*}
		\begin{aligned}
			&\left\|\nabla f(Z^{k+1})\!-\! \nabla f(Z^{k+\frac{1}{2}})\!+\! \nabla_{Z}H(X^{k+1}\!,\!Y^{k+1}\!,\!Z^{k+1})\!-\! \nabla_{Z}H(X^{k+1}\!,\!Y^{k}\!,\!Z^{k+\frac{1}{2}})\right\|_F \\
			\overset{\text{(a6)}}{\le}& M\left\|Y^{k+1}-Y^k\right\|_F+\left(M+L_f\right)\left\|Z^{k+1}-Z^{k+\frac12}\right\|_F.
		\end{aligned}
	\end{equation*}
	The RHS of \eqref{lem33-5}
	\begin{equation*}
		\begin{aligned}
			& \beta \left \| \mathcal{C}^*\left [\mathcal{B}\left( Y^{k+1}-Y^k\right) +\mathcal{C}\left( Z^{k+1}-Z^{k+\frac12}\right) \right ] \right \|_F \\
			\overset{\text{(a2)}}{\ge}& \frac{\beta}{\sqrt{\mu_4}}
			\left\|\mathcal{B}\left( Y^{k+1}-Y^k\right) +\mathcal{C}\left( Z^{k+1}-Z^{k+\frac12}\right)  \right\|_F \\
			\ge \; & \frac{\beta}{\sqrt{\mu_4}}
			\left | \left\| \mathcal{B}\left( Y^{k+1}-Y^k\right) \right \|_F -  \left\|  \mathcal{C}\left( Z^{k+1}-Z^{k+\frac12}\right) \right \|_F \right |.
		\end{aligned}
	\end{equation*}
	Combining the above three expressions, we obtain
	\begin{equation}\label{lem33-8}
		\begin{aligned}
			&M\left\|Y^{k+1}-Y^k\right\|_F+\left(M+L_f\right)\left\|Z^{k+1}-Z^{k+\frac12}\right\|_F \\
			\ge \;&\frac{\beta}{\sqrt{\mu_4}}
			\left |  \left\| \mathcal{B}\left( Y^{k+1}-Y^k\right)  \right \|_F -  \left\|  \mathcal{C}\left( Z^{k+1}-Z^{k+\frac12}\right)  \right\|_F \right |.
		\end{aligned}
	\end{equation}
	If $\left\| \mathcal{B}\left( Y^{k+1}-Y^k\right) \right\|_F < \left\|  \mathcal{C}\left( Z^{k+1}-Z^{k+\frac12}\right) \right\|_F$, then \eqref{lem33-8} tells that
	\begin{equation*}
		\begin{aligned}
			&\left ( \beta\sqrt{\frac{\mu_3}{\mu_4}}-M-L_f \right )\left\|Z^{k+1}-Z^{k+\frac12}\right\|_F \\
			\overset{\text{(a2)}}{\le}\hspace{0.4em}& \frac{\beta}{\sqrt{\mu_4}}\left\|\mathcal{C}\left( Z^{k+1}-Z^{k+\frac12}\right) \right\|_F-\left(M+L_f\right)\left\|Z^{k+1}-Z^{k+\frac12}\right\|_F \\
			\overset{\eqref{lem33-8}}{\le} & \left (\frac{\beta}{\sqrt{\mu_4}} \|\mathcal{B}\| +M  \right )\left\|Y^{k+1}-Y^k\right \|_F. 
		\end{aligned}
	\end{equation*}
	Hence, we have
	\begin{equation}\label{equ:Z}
		\left\|Z^{k+1}-Z^{k+\frac12}\right\|_F\le \frac{\beta\|\mathcal{B}\|+\sqrt{\mu_4}M}{\beta\sqrt{\mu_3}-(M+L_f)\sqrt{\mu_4}}\left\| Y^{k+1}-Y^k \right\|_F
	\end{equation}
	since $\frac{\beta\|\mathcal{B}\|+\sqrt{\mu_4}M}{\beta\sqrt{\mu_3}-(M+L_f)\sqrt{\mu_4}}\overset{\text{(a7)}}{\ge} \frac{\|\mathcal{B}\|}{\sqrt{\mu_3}}$.
	
	At last, combining \eqref{equ:lambda}, \eqref{lem33-4}, and \eqref{equ:Z}, we get the desired boundness
	\begin{align*}
		\mathrm{dist}\left( 0,\partial \mathcal{L}_{\beta}(W^{k+1})\right)  
		& \le \vartheta_1 \left\| Y^{k+1}-Y^k \right\|_F + \vartheta_2 \left\| Z^{k+1}-Z^k \right\|_F \\
		& \le \vartheta \left(\left\| Y^{k+1}-Y^k \right\|_F+\left\| Z^{k+1}-Z^k \right\|_F\right),
	\end{align*}
	where 
	\begin{align*}
		&\vartheta_1:=\left(\|\mathcal{A}^*\|+\|\mathcal{B}^*\|+\|\mathcal{C}^*\|+\frac{1}{\beta}\right)\frac{M+L_f}{\sqrt{\mu_4}}+ \frac{(\beta\|\mathcal{B}^*\|\|\mathcal{C}\|+L_1)(\beta\|\mathcal{B}\|+\sqrt{\mu_4}M)}{\beta\sqrt{\mu_3}-(M+L_f)\sqrt{\mu_4}} \\
		& \hspace{2.6em} +\beta\|\mathcal{A}^*\|\|\mathcal{B}\| +N>0, \\
		&\vartheta_2:=\left(\|\mathcal{A}^*\|+\|\mathcal{B}^*\|+\|\mathcal{C}^*\|+\frac{1}{\beta}\right)\frac{M+L_f}{\sqrt{\mu_4}}+\beta\|\mathcal{A}^*\|\|\mathcal{C}\| +N>0, \\
		&\vartheta:=\max \{ \vartheta_1, \vartheta_2\}.
	\end{align*}
\end{proof}

\begin{lemma}[Subsequence convergence]\label{lem:subconverge}
	Suppose $\left\{W^k\right\}_{k\in \mathbb{N}}$ generated by ADMMn is bounded and Assumption \ref{ass:converge} holds. Denote the set of limit points of $\left\{W^k\right\}$ as $\mathcal{S}(W^0)$. Then, we have
	\begin{enumerate}
		\item[(1)] $\mathcal{S}(W^0)$ is nonempty, compact and connected. Moreover,
		\begin{equation}\label{lem34-1}
			\mathrm{dist}\left( W^k,\mathcal{S}(W^0)\right) \to 0 \; \text{ as } k\to \infty.
		\end{equation}
		\item[(2)] $\mathcal{S}(W^0)\subseteq \mathrm{crit}\mathcal{L}_{\beta}$.
		\item[(3)] $\mathcal{L}_{\beta}$ equals to $\lim_{k\to \infty}\mathcal{L}_{\beta}(W^k)$ on $\mathcal{S}(W^0)$.
	\end{enumerate}
\end{lemma}

\begin{proof}
	We prove the above statements one by one.
	
	(1) From the definition of the set of limit points, it is trivial and can be found in \cite[Lemma 5]{limitset}.
	
	(2) For any $W^{\star}\in \mathcal{S}(W^0)$, there exists a subsequence $\left\{W^{k_j}\right\}$ that converges to $W^{\star}$. From Lemma \ref{lem:square}, we know $\left\|W^{k+1}-W^k\right\|_F^2\to 0$, which means that $\left\{W^{k_j+1}\right\}$ also converges to $W^{\star}$.
	Thus, from
	\begin{equation*}
		\mathcal{L}_{\beta}\left(X^{k+1},Y^k,Z^k,\Lambda^k\right)\overset{\eqref{equ:admmscheme1}}{\le} \mathcal{L}_{\beta}\left(X^{\star},Y^k,Z^k,\Lambda^k\right),
	\end{equation*}
	we conclude that
	\begin{align*}
		\limsup_{j\to \infty}\mathcal{L}_{\beta}\left(W^{k_j+1}\right)&=\limsup_{j\to \infty}\mathcal{L}_{\beta}\left(X^{k_j+1},Y^{k_j},Z^{k_j},\Lambda^{k_j}\right)\\
		&\le \limsup_{j\to \infty}\mathcal{L}_{\beta}\left(X^{\star},Y^{k_j},Z^{k_j},\Lambda^{k_j}\right)=\mathcal{L}_{\beta}(W^{\star}),
	\end{align*}
	where both equalities are due to the continuity of $\mathcal{L}_{\beta}$ w.r.t. $(Y,Z,\Lambda)$. 
	Since $\mathcal{L}_{\beta}$ is lower semicontinuous, we also have
	\begin{equation*}
		\liminf_{j\to \infty}\mathcal{L}_{\beta}\left(W^{k_j+1}\right)\ge \mathcal{L}_{\beta}(W^{\star}).
	\end{equation*}
	The above two inequalities tell that
	\begin{equation}\label{lem34-4}
		\lim_{j\to \infty}\mathcal{L}_{\beta}\left(W^{k_j+1}\right)= \mathcal{L}_{\beta}(W^{\star}).
	\end{equation}
	
	Lemma \ref{lem:square} implies $\left\|Y^{k_j+1}-Y^{k_j}\right\|_F\to 0$ and  $\left\|Z^{k_j+1}-Z^{k_j}\right\|_F\to 0$. Then it follows that $\mathrm{dist}\left( 0,\partial \mathcal{L}_{\beta}(W^{k_j+1})\right) \to 0$ by Lemma \ref{lem:boundsub}. Since $\partial \Phi$, $\partial \Psi$ are closed, $\nabla f$, $\nabla_X H$, $\nabla_Y H$ and $\nabla_Z H$ are continuous, taking limit along the subsequence $\left\{W^{k_j+1}\right\}$, we obtain $(0,0,0,0)\in \partial \mathcal{L}_{\beta}(W^{\star})$. Hence $W^{\star} \in \mathrm{crit}\mathcal{L}_{\beta}$, which proves $\mathcal{S}(W^0)\subseteq \mathrm{crit}\mathcal{L}_{\beta}$.
	
	(3) Lemma \ref{lem:descent} implying that $\left\{\mathcal{L}_{\beta}(W^k)\right\}_{k\in \mathbb{N}}$ is monotonically decreasing, combining \eqref{lem34-4}, we immediately have 
	\begin{equation}\label{lem34-5}
		\lim_{k\to \infty}\mathcal{L}_{\beta}(W^{k})= \mathcal{L}_{\beta}(W^{\star}).
	\end{equation}
	Hence $\mathcal{L}_{\beta}$ equals to the constant $\lim_{k\to \infty}\mathcal{L}_{\beta}(W^{k})$ on $\mathcal{S}(W^0)$.
\end{proof}

Now, we are ready to give the main convergence result.

\begin{theorem}[Global convergence]\label{thm:converge}
	Suppose $\mathcal{L}_{\beta}$ is a KL function and Assumption \ref{ass:converge} holds. Let $\left\{W^k\right\}_{k\in \mathbb{N}}$ be the bounded sequence generated by ADMMn. Then, we have 
	\begin{equation}\label{equ:finitelength}
		\sum_{k=0}^{\infty}\left\|W^{k+1}-W^k\right\|_F<+\infty.
	\end{equation}
	As a consequence, $\left\{W^k\right\}$ globally converges to $W^{\star}\in \mathrm{crit}\mathcal{L}_{\beta}$.
\end{theorem}

\begin{proof}
	Lemma \ref{lem:subconverge} implies that \eqref{lem34-5} holds. Hence we just discuss two cases.
	
	(1) Case 1: there exists some $k_0\in \mathbb{N}$ such that $\mathcal{L}_{\beta}(W^{k_0})=\mathcal{L}_{\beta}(W^{\star})$. 
	
	From Lemma \ref{lem:descent}, we know for any $k>k_0$, 
	\begin{align*}
		&c_1\left\|Y^{k+1}-Y^k\right\|_F^2+c_2\left\|Z^{k+1}-Z^k\right\|_F^2 \\
		\le \;& \mathcal{L}_\beta(W^k)-\mathcal{L}_\beta(W^{k+1})\le \mathcal{L}_\beta(W^{k_0})-\mathcal{L}_\beta(W^{\star})=0,
	\end{align*}
	hence $Y^{k+1}=Y^k$ and $Z^{k+1}=Z^k$. Together with \eqref{equ:lambda} we immediately get $\Lambda^{k+1}=\Lambda^k$. Combining \eqref{equ:X} we know $X^{k+1}=X^k$ for any $k>k_0+1$. Therefore, $W^{k+1}=W^k$ for any $k>k_0+1$, which proves \eqref{equ:finitelength}.
	
	(2) Case 2: for any $k\in \mathbb{N}$, $\mathcal{L}_{\beta}(W^{k})>\mathcal{L}_{\beta}(W^{\star})$. 
	
	From \eqref{lem34-1} and \eqref{lem34-5}, we know for any $\eta$, $\sigma>0$, there exists $k_1\in \mathbb{N}$ such that $\mathcal{L}_{\beta}(W^k)<\mathcal{L}_{\beta}(W^{\star})+\eta$ and $\mathrm{dist}\left( W^k,\mathcal{S}(W^0)\right) <\sigma$ for any $k>k_1$.
	Besides, Lemma \ref{lem:subconverge} implies that $\mathcal{S}(W^0)$ is a nonempty compact set, and $\mathcal{L}_{\beta}$ is a constant on it. Let $\Omega=\mathcal{S}(W^0)$ in Proposition \ref{pro:UKL}, we know for any $k>k_1$,
	\begin{equation}\label{thm1-1}
		\varphi^{\prime}\left(\mathcal{L}_{\beta}(W^k)-\mathcal{L}_{\beta}(W^{\star})\right)\cdot \mathrm{dist}\left(0,\partial \mathcal{L}_{\beta}(W^k)\right)\geq 1.
	\end{equation} 
	For convenience, denote $\Delta_{k,k+1}:=\varphi(\mathcal{L}_{\beta}(W^k)-\mathcal{L}_{\beta}(W^{\star}))-\varphi(\mathcal{L}_{\beta}(W^{k+1})-\mathcal{L}_{\beta}(W^{\star}))$. Since $\varphi$ is concave, we have
	\begin{equation*}
		\Delta_{k,k+1}\ge \varphi^{\prime}\left( \mathcal{L}_{\beta}(W^k)-\mathcal{L}_{\beta}(W^{\star})\right)  \cdot \left( \mathcal{L}_{\beta}(W^k)-\mathcal{L}_{\beta}(W^{k+1})\right) .
	\end{equation*}
	Combining $\varphi^{\prime}>0$ on $(0,\eta)$, it follows that
	\begin{equation}\label{thm1-2}
		\mathcal{L}_{\beta}(W^k)-\mathcal{L}_{\beta}(W^{k+1})\le \frac{\Delta_{k,k+1}}{\varphi^{\prime}\left( \mathcal{L}_{\beta}(W^k)-\mathcal{L}_{\beta}(W^{\star})\right) }.
	\end{equation}
	The LHS of \eqref{thm1-2}
	\begin{align*}
		\mathcal{L}_{\beta}(W^k)-\mathcal{L}_{\beta}(W^{k+1})
		\overset{\eqref{equ:descent}}{\ge}& c\left(\left\|Y^{k+1}-Y^k\right\|_F^2+\left\|Z^{k+1}-Z^k\right\|_F^2\right) \\
		\ge \hspace{0.5em} & \frac{c}{2}\left(\left\|Y^{k+1}-Y^k\right\|_F+\left\|Z^{k+1}-Z^k\right\|_F\right)^2,
	\end{align*}
	where $c:=\min\{c_1,c_2\}>0$ in (\ref{equ:descent}). \\
	The RHS of \eqref{thm1-2}
	\begin{equation}\label{thm1-4}
		\begin{aligned}
			\frac{\Delta_{k,k+1}}{\varphi^{\prime}\left( \mathcal{L}_{\beta}(W^k)-\mathcal{L}_{\beta}(W^{\star})\right) } 
			&\overset{\eqref{thm1-1}}{\le} \Delta_{k,k+1}\cdot \mathrm{dist}\left( 0,\partial \mathcal{L}_{\beta}(W^k)\right)  \\
			&\overset{\eqref{equ:boundsub}}{\le} \vartheta\left(\left\|Y^k-Y^{k-1}\right\|_F+\left\|Z^k-Z^{k-1}\right\|_F\right)\cdot \Delta_{k,k+1}.
		\end{aligned}
	\end{equation}
	The above three inequalities imply that for any $k>k_1$, 
	\begin{align*}
		&\left(\left\|Y^{k+1}-Y^k\right\|_F+\left\|Z^{k+1}-Z^k\right\|_F\right)^2 \\ \le \;& \frac{2\vartheta}{c}\left(\left\|Y^k-Y^{k-1}\right\|_F+\left\|Z^k-Z^{k-1}\right\|_F\right)\cdot \Delta_{k,k+1}.
	\end{align*}
	Thus,
	\begin{equation*}
		\begin{aligned}
			& \left\|Y^{k+1}-Y^k\right\|_F+\left\|Z^{k+1}-Z^k\right\|_F \\
			\le \;&\sqrt{\frac{2\vartheta}{c}\Delta_{k,k+1}}\cdot \sqrt{\left\|Y^k-Y^{k-1}\right\|_F+\left\|Z^k-Z^{k-1}\right\|_F} \\
			\le\;& \frac{\vartheta}{c}\Delta_{k,k+1}+\frac12 \left(\left\|Y^k-Y^{k-1}\right\|_F+\left\|Z^k-Z^{k-1}\right\|_F\right).
		\end{aligned}
	\end{equation*}
	Summing up from $k=k_1+1$ to $K$, then letting $K\to \infty$, we conclude
	\begin{equation}\label{thm1-5}
		\begin{aligned}
			&\sum_{k=k_1+1}^{\infty}\left(\left\|Y^{k+1}-Y^k\right\|_F+\left\|Z^{k+1}-Z^k\right\|_F\right)\\
			\le~& \left\|Y^{k_1+1}-Y^{k_1}\right\|_F+\left\|Z^{k_1+1}-Z^{k_1}\right\|_F + \frac{2\vartheta}{c}\varphi\left( \mathcal{L}_{\beta}(W^{k_1+1})-\mathcal{L}_{\beta}(W^{\star})\right) <+\infty,
		\end{aligned}
	\end{equation}
	which means
	\begin{equation*}
		\sum_{k=0}^{\infty}\left\|Y^{k+1}-Y^k\right\|_F<+\infty, \quad  \sum_{k=0}^{\infty}\left\|Z^{k+1}-Z^k\right\|_F<+\infty.
	\end{equation*}
	Besides,
	\begin{equation*}
		\sum_{k=0}^{\infty}\left\|\Lambda^{k+1}-\Lambda^k\right\|_F \overset{\eqref{equ:lambda}}{\le} \frac{M+L_f}{\sqrt{\mu_4}} \sum_{k=0}^{\infty}\left(\left\|Y^{k+1}-Y^k\right\|_F+\left\|Z^{k+1}-Z^k\right\|_F\right)<+\infty,
	\end{equation*}
	and
	\begin{align*}
		& \sum_{k=0}^{\infty}\left\|X^{k+1}-X^k\right\|_F \\
		\overset{\eqref{equ:X}}{\le}& \left\|X^{1}-X^0\right\|_F+
		\frac{1}{\beta \sqrt{\mu_1}}\sum_{k=1}^{\infty}\left(\left\|\Lambda^{k}-\Lambda^{k-1}\right  \|_F+ \left\|\Lambda^{k+1}-\Lambda^{k}\right  \|_F\right) \\
		~&+\frac{1}{\sqrt{\mu_1}}\sum_{k=1}^{\infty}\left(\|\mathcal{B}\|\left\|Y^{k+1}-Y^{k} \right \|_F+  \|\mathcal{C}\| \left\|Z^{k+1}-Z^{k}\right \|_F\right)<+\infty.
	\end{align*}
	Therefore, we prove $\sum_{k=0}^{\infty}\left\|W^{k+1}-W^k\right\|_F<+\infty$.
	
	Above all, no matter which cases, it holds that $\sum_{k=0}^{\infty}\left\|W^{k+1}-W^k\right\|_F<+\infty$, which means $\left\{W^k\right\}_{k\in \mathbb{N}}$ is a Cauchy sequence. Hence $\{W^k\}$ is globally convergent. By Lemma \ref{lem:subconverge}, we know that $W^k \to W^{\star}\in \mathrm{crit}\mathcal{L}_{\beta}$.
	This completes the proof.
\end{proof}

\begin{theorem}[Convergence rate]\label{thm:rate}
	Let $\left\{W^k\right\}_{k\in \mathbb{N}}$ be the sequence generated by ADMMn which converges to $W^{\star}$. Assume $\mathcal{L}_{\beta}$ has the KL property at $W^{\star}$ with $\varphi(x)=tx^{1-\theta}$, $\theta\in [0,1)$, $t>0$. Then the convergence rate of ADMMn is determined by $\theta$ as follows:
	\begin{enumerate}
		\item[(1)] If $\theta=0$, then $\left\{W^k\right\}$ converges to $W^{\star}$ in a finite number of steps.
		\item[(2)] If $\theta\in \left( 0,\frac12\right] $, then $\exists \; \Bar{c}_1>0$, $\tau\in [0,1)$ such that $\left\|W^k-W^{\star}\right\|_F\le \Bar{c}_1\tau^k$, i.e., the convergence rate is linear.
		\item[(3)] If $\theta\in \left( \frac12,1\right) $, then $\exists \; \Bar{c}_2>0$ such that $\left\|W^k-W^{\star}\right\|_F\le\Bar{c}_2k^{\frac{\theta-1}{2\theta-1}}$, i.e., the convergence rate is sublinear.
	\end{enumerate}  
\end{theorem}

\begin{proof}
	We prove the theorem based on \cite[Theorem 2]{rate}.
	
	(1) If $\theta=0$, then $\varphi(x)=tx$, $\varphi^{\prime}(x)=t$. Assuming that $\left\{W^k\right\}$ does not converge to $W^{\star}$ in a finite number of steps, then \eqref{thm1-1} tells $t \cdot \mathrm{dist}\left( 0,\partial \mathcal{L}_{\beta}(W^k)\right) \ge 1$ for any $k>k_1$, which is contradict to Lemma \ref{lem:boundsub}.
	
	Now suppose $\theta>0$ and denote $\Delta_k:=\sum_{i=k}^{\infty}\left(\left\|Y^{i+1}-Y^i\right\|_F+\left\|Z^{i+1}-Z^i\right\|_F\right)$ for $k\in \mathbb{N}$. 
	\eqref{thm1-5} implies 
	\begin{equation}\label{thm2-1}
		\Delta_{k_1+1}\le \Delta_{k_1}-\Delta_{k_1+1}+\frac{2\vartheta}{c}\varphi\left(\mathcal{L}_{\beta}(W^{k_1+1})-\mathcal{L}_{\beta}(W^{\star})\right).
	\end{equation}
	Since $\varphi(x)=tx^{1-\theta}$, \eqref{thm1-1} tells
	\begin{equation*}
		t(1-\theta)\cdot \mathrm{dist}\left( 0,\partial \mathcal{L}_{\beta}(W^{k_1+1})\right)  \ge \left(\mathcal{L}_{\beta}(W^{k_1+1})-\mathcal{L}_{\beta}(W^{\star})\right)^{\theta}.
	\end{equation*}
	Besides, \eqref{equ:boundsub} implies
	\begin{equation*}
		\mathrm{dist}\left( 0,\partial \mathcal{L}_{\beta}(W^{k_1+1})\right) \le \vartheta \left(\left\|Y^{k_1+1}\!-\! Y^{k_1}\right\|_F+\left\|Z^{k_1+1}\!-\! Z^{k_1}\right\|_F\right)=\vartheta(\Delta_{k_1}-\Delta_{k_1+1}).
	\end{equation*}
	Combining the above two inequalities and $\varphi(x)=tx^{1-\theta}$, we have
	\begin{align*}
		& \varphi\left(\mathcal{L}_{\beta}(W^{k_1+1})-\mathcal{L}_{\beta}(W^{\star})\right) 
		= t\left(\mathcal{L}_{\beta}(W^{k_1+1})-\mathcal{L}_{\beta}(W^{\star})\right)^{1-\theta} \\
		\le \; & t^{\frac{1}{\theta}}(1-\theta)^{\frac{1-\theta}{\theta}}\cdot \left(\mathrm{dist}\left( 0,\partial \mathcal{L}_{\beta}(W^{k_1+1})\right) \right)^{\frac{1-\theta}{\theta}} 
		\le \gamma (\Delta_{k_1}-\Delta_{k_1+1})^{\frac{1-\theta}{\theta}},
	\end{align*}
	where $\gamma:=t^{\frac{1}{\theta}}\left [ (1-\theta)\vartheta \right ]^{\frac{1-\theta}{\theta}}>0$.
	Substituting it into \eqref{thm2-1} yields
	\begin{equation*}
		\Delta_{k_1+1}\le \Delta_{k_1}-\Delta_{k_1+1}+\frac{2\vartheta}{c}\gamma (\Delta_{k_1}-\Delta_{k_1+1})^{\frac{1-\theta}{\theta}}.
	\end{equation*}
	
	The inequality of this form has been studied in \cite[Theorem 2]{rate}:
	
	(2) If $\theta\in (0,1/2]$, then $\exists \; \Bar{c}>0$, $\tau\in [0,1)$ such that (cf. \cite[Theorem 2]{rate})
	\begin{equation*}
		\left\|Y^{k}-Y^{\star}\right\|_F+\left\|Z^{k}-Z^{\star}\right\|_F \le \Delta_k \le \Bar{c}\tau^k.
	\end{equation*}
	Combining \eqref{equ:1optadmm3} and \eqref{equ:2optadmm6}, similar to the process of estimating \eqref{equ:lambda}, we have
	\begin{equation*}
		\left\| \Lambda^{k}-\Lambda^{\star}\right \|_F \le \frac{M+L_f}{\sqrt{\mu_4}}\left(\left\| Y^{k}-Y^{\star}\right\|_F+\left\|Z^{k}-Z^{\star} \right\|_F\right).
	\end{equation*}
	Combining \eqref{equ:1optadmm4} and \eqref{equ:2optadmm5}, analogous to the estimation of \eqref{equ:X}, we get
	\begin{align*}
		\left\|X^{k}-X^{\star}\right\|_F\le \;& \frac{1}{\beta \sqrt{\mu_1}}\left(\left\|\Lambda^{k}-\Lambda^{\star} \right \|_F+\left \|\Lambda^{k-1}-\Lambda^{\star} \right \|_F \right. \\
		& \left. + \beta \|\mathcal{B}\|\left\|Y^{k}-Y^{\star} \right \|_F+ \beta \|\mathcal{C}\| \left\|Z^{k}-Z^{\star}\right \|_F\right).
	\end{align*}
	The above three inequalities tell that 
	\begin{align*}
		&\left\|W^{k}-W^{\star}\right\|_F 
		\le \left\|X^{k}-X^{\star}\right\|_F+\left\|Y^{k}-Y^{\star}\right\|_F+\left\|Z^{k}-Z^{\star}\right\|_F+\left\| \Lambda^{k}-\Lambda^{\star} \right\|_F \\
		\le \; & \gamma_1\left(\left\| Y^{k-1}-Y^{\star}\right\|_F+\left\|Z^{k-1}-Z^{\star} \right\|_F\right) +\gamma_1\left(\left\| Y^{k}-Y^{\star}\right\|_F+\left\|Z^{k}-Z^{\star}\right \|_F\right) \\
		\le \; & \gamma_1\Bar{c}\left(\tau^{k-1}+\tau^{k}\right) 
		= \Bar{c}_1\tau^{k},
	\end{align*}
	where 
	\begin{align*}
		&\gamma_1:=\left (\frac{1}{\beta \sqrt{\mu_1}}+1  \right )\frac{M+L_f}{\sqrt{\mu_4}}+\frac{1}{\sqrt{\mu_1}}\max\left\lbrace  \|\mathcal{B}\|, \|\mathcal{C}\| \right\rbrace >0, \\
		& \Bar{c}_1:=\gamma_1\Bar{c}\left (1+\frac{1}{\tau}  \right )>0.
	\end{align*}
	
	(3) If $\theta\in (1/2,1)$, then $\exists \; \Bar{c}_2>0$, such that (cf. \cite[Theorem 2]{rate})
	\begin{equation*}
		\left\|Y^{k}-Y^{\star}\right\|_F+\left\|Z^{k}-Z^{\star}\right\|_F \le \Delta_k \le \Bar{c}k^{\frac{\theta-1}{2\theta-1}}.
	\end{equation*}
	Similar to (2), we have
	\begin{align*}
		\left\|W^{k}-W^{\star}\right\|_F 
		\le \;&\gamma_1\left(\left\| Y^{k-1}-Y^{\star}\right\|_F+\left\|Z^{k-1}-Z^{\star} \right\|_F\right) \\
		& +\gamma_1\left(\left\| Y^{k}-Y^{\star}\right\|_F+\left\|Z^{k}-Z^{\star}\right \|_F\right) \\
		\le\; & \gamma_1\Bar{c} \left ( (k-1)^{\frac{\theta-1}{2\theta-1}}+k^{\frac{\theta-1}{2\theta-1}} \right ) 
		\le \Bar{c}_2 k^{\frac{\theta-1}{2\theta-1}},
	\end{align*}
	where $\Bar{c}_2:=2\gamma_1 \Bar{c}>0$.
	
	Above all, we complete the proof.
\end{proof}

The following proposition provides a sufficient condition to obtain a bounded sequence $\left\{W^k\right\}$ generated by ADMMn. Note that there are other conditions to bound $\left\{W^k\right\}$. We can also use similar ways to Proposition \ref{pro:boundW} to check if the sequence $\left\{W^k\right\}$ generated by ADMMn is bounded for any particular problems.

\begin{prop}[Boundness of $\left\{W^k\right\}$]\label{pro:boundW}
	Suppose $\left\{W^k\right\}$ is the sequence generated by ADMMn. If
	\begin{enumerate}
		\item[(1)] $\inf_{X}\Phi(X)>-\infty$ and $\inf_{Y}\Psi(Y)>-\infty$;
		\item[(2)] f(Z) is coercive, i.e., $\liminf_{\|Z\|_F\to \infty}f(Z)=+\infty$, and
		\begin{equation*}
			\Bar{f}:=\inf_Z \left\{ f(Z)-\frac{2}{\mu_4\hat{\beta}}\|\nabla f(Z)\|_F^2 \right\}>-\infty;
		\end{equation*}
		\item[(3)] H(X,Y,Z) is coercive to $X$ and $Y$, i.e., $\liminf_{\|X\|_F+\|Y\|_F \to \infty}H(X,Y,Z)=+\infty$, 
		and 
		\begin{equation*}
			\Bar{H}:=\inf_{X,Y,Z} \left\{ H(X,Y,Z)-\frac{2}{\mu_4\hat{\beta}}\|\nabla_Z H(X,Y,Z)\|_F^2 \right\}>-\infty.
		\end{equation*}
	\end{enumerate}
	Then, the sequence $\left\{W^k\right\}$ is bounded.
\end{prop}

\begin{proof}
	From \eqref{equ:2optadmm6}, we have
	\begin{equation}\label{equ:pro31-1}
		\begin{aligned}
			\left\|\Lambda^k\right\|_F^2\overset{\text{(a2)}}{\le} & \frac{1}{\mu_4}\left\|\mathcal{C}^*(\Lambda^k)\right\|_F^2\overset{\eqref{equ:2optadmm6}}{=}\frac{1}{\mu_4}\left\|\nabla f(Z^k)+\nabla_Z H(X^k,Y^k,Z^k)\right\|_F^2\\
			\le \;& \frac{2}{\mu_4}\left(\left\|\nabla f(Z^k)\right\|_F^2+\left\|\nabla_Z H(X^k,Y^k,Z^k)\right\|_F^2\right).
		\end{aligned}
	\end{equation}
	It follows from Lemma \ref{lem:descent} that
	\begin{align*}
		& \mathcal{L}_{\beta}(W^1) \ge \mathcal{L}_{\beta}(W^k)\\
		=\;&\Phi(X^k)+\Psi(Y^k)+f(Z^k)+H(X^k,Y^k,Z^k)\\
		&-\left\langle \Lambda^k,\mathcal{A}(X^k)+\mathcal{B}(Y^k)+\mathcal{C}(Z^k)-D\right\rangle+\frac{\beta}{2}\left\|\mathcal{A}(X^k)+\mathcal{B}(Y^k)+\mathcal{C}(Z^k)-D\right\|_F^2 \\
		=\;&\Phi(X^k)+\Psi(Y^k)+f(Z^k)+H(X^k,Y^k,Z^k)\\
		&+\frac{\beta}{2}\left\|\mathcal{A}(X^k)+\mathcal{B}(Y^k)+\mathcal{C}(Z^k)-D-\frac{\Lambda^k}{\beta}\right\|_F^2-\frac{1}{2\beta}\left\|\Lambda^k\right\|_F^2  \\
		=\;&\Phi(X^k)+\Psi(Y^k)+\frac12 f(Z^k)+\frac12 \left ( f(Z^k)-\frac{2}{\mu_4\hat{\beta}}\left\|\nabla f(Z^k)\right\|_F^2 \right )\\
		&+\frac12 H(X^k,Y^k,Z^k)+\frac12 \left ( H(X^k,Y^k,Z^k)-\frac{2}{\mu_4\hat{\beta}}\left\|\nabla_Z H(X^k,Y^k,Z^k)\right\|_F^2 \right )\\
		&+\frac{\beta}{2}\left\|\mathcal{A}(X^k)+\mathcal{B}(Y^k)+\mathcal{C}(Z^k)-D-\frac{\Lambda^k}{\beta}\right\|_F^2 \\
		&+\frac{1}{\mu_4\hat{\beta}}\left(\left\|\nabla f(Z^k)\right\|_F^2+\left\|\nabla_Z H(X^k,Y^k,Z^k)\right\|_F^2\right)
		-\frac{1}{2\beta}\left\|\Lambda^k\right\|_F^2 \\
		\overset{\eqref{equ:pro31-1}}{\ge}& \Phi(X^k)+\Psi(Y^k)+\frac12 f(Z^k) +\frac12 \Bar{f}+\frac12 H(X^k,Y^k,Z^k)+\frac12 \Bar{H} \\
		&+\frac{\beta}{2}\left\|\mathcal{A}(X^k)+\mathcal{B}(Y^k)+\mathcal{C}(Z^k)-D-\frac{\Lambda^k}{\beta}\right\|_F^2+\frac12 \left (\frac{1}{\hat{\beta}}- \frac{1}{\beta} \right ) \left\|\Lambda^k\right\|_F^2.
	\end{align*}
	Due to conditions (1), (2), (3), and $\beta>\hat{\beta}$ in Assumption \ref{ass:converge} (a7), we know that $\left\{X^k\right\}$, $\left\{Y^k\right\}$, $\left\{Z^k\right\}$, $\left\{\Lambda^k\right\}$ and $\left\{\mathcal{A}(X^k)+\mathcal{B}(Y^k)+\mathcal{C}(Z^k)-D-\frac{\Lambda^k}{\beta}\right\}$ are all bounded. Thus $\left\{W^k\right\}$ is bounded.
	
	In particular, if $H(X,Y,Z)$ is a function independent of $X$ or $Y$, e.g., $H(X , Y , Z)= H(Y,Z)$, the above argument can only tell that $\left\{\mathcal{A}(X^k)\!+\!\mathcal{B}(Y^k)\!+\!\mathcal{C}(Z^k)\!-\! D \!-\!\frac{\Lambda^k}{\beta}\right\}$, $\left\{Y^k\right\}$, $\left\{Z^k\right\}$ and $\left\{\Lambda^k\right\}$ are bounded. In this case,
	denote $E^k:=\mathcal{A}(X^k)+\mathcal{B}(Y^k)+\mathcal{C}(Z^k)-D-\frac{\Lambda^k}{\beta}$. From
	\begin{align*}
		\left\|X^k\right\|_F^2 \overset{\text{(a2)}}{\le}& \frac{1}{\mu_1}\left\|\mathcal{A}(X^k)\right\|_F^2=\frac{1}{\mu_1}\left\|E^k -\mathcal{B}(Y^k)-\mathcal{C}(Z^k)+D+\frac{\Lambda^k}{\beta}\right\|_F^2\\
		\le \;&\frac{5}{\mu_1}\left(\left\|E^k\right\|_F^2+\left\|\mathcal{B}(Y^k)\right\|_F^2+\left\|\mathcal{C}(Z^k)\right\|_F^2+\left\|D\right\|_F^2+\left\|\frac{\Lambda^k}{\beta}\right\|_F^2\right) <+\infty,
	\end{align*}
	we also get $\left\{X^k\right\}$ bounded. Therefore, the sequence $\left\{W^k\right\}$ is bounded.
\end{proof}

\begin{remark}
	In practice, we may also encounter situations where $f(Z)$ is not coercive. In this case, using analogous proofs in Proposition \ref{pro:boundW}, we are able to obtain another sufficient condition to bound $\left\{W^k\right\}$ with a stronger assumption on $H(X,Y,Z)$, that is $\liminf_{\|X\|_F+\|Y\|_F+\|Z\|_F \to \infty}H(X,Y,Z)=+\infty$.
\end{remark}

At the end of this section, we discuss two scenarios in applications of ADMMn and the corresponding countermeasures. First, when $\Psi \equiv 0$, \eqref{equ:opt} degenerates into a 2-block nonconvex problem. At this time, ADMMn also becomes the classic 2-block ADMM algorithm. We may still want to apply ADMMn with $Z$ being updated twice at each iteration, which can be achieved by swapping the update order of $X$ and $Y$ in \eqref{equ:admmscheme1}-\eqref{equ:admmscheme5}. However, a natural question arises: is ADMMn still convergent with the update order of $X$ and $Y$ being swapped? The answer is yes. Here is the proposition.

\begin{prop}\label{pro:swap}
	Swapping the update order of $X$ and $Y$ in \eqref{equ:admmscheme1}-\eqref{equ:admmscheme5}, ADMMn still has the global convergence under the following simple modifications to (a3), (a4), and (a5) in Assumption \ref{ass:converge}.
	\begin{enumerate}
		\item[(a3)] There exists $N(Y)$ such that $\|\nabla_{Y}H(X_{1},Y,Z_{1})-\nabla_{Y}H(X_{2},Y,Z_{2})\|_F\leq N(Y)(\|X_{1}-X_{2}\|_F+\|Z_{1}-Z_{2}\|_F)$ for any fixed $Y$;
		\item[(a4)] $\nabla_X H(X,Y,\cdot)$ is $L_1(X,Y)$-Lipschitz for any fixed $X$ and $Y$.  $\nabla_Y H(X,\cdot,Z)$ is $L_2(X,Z)$-Lipschitz for any fixed $X$ and $Z$. $\nabla_Z H(X,Y,\cdot)$ is $L_3(X,Y)$-Lipschitz for any fixed $X$ and $Y$; 
		\item[(a5)] There exist $L_1$, $L_2$, $L_3>0$, $N>0$ such that $\sup_{k\in \mathbb{N}} L_1\left(X^{k+1},Y^{k+1}\right)\le L_1$, $\sup_{k\in \mathbb{N}} L_2\left(X^{k},Z^{k}\right)\le L_2$, $\sup_{k\in \mathbb{N}} \left\{L_3\left(X^{k},Y^{k}\right),L_3\left(X^{k},Y^{k+1}\right)\right\}  \le L_3$ and $\sup_{k\in \mathbb{N}} N\left(Y^{k+1}\right)\le N$.
	\end{enumerate}
\end{prop}

\begin{proof}
	The proof is highly similar to the convergence analysis of ADMMn \eqref{equ:admmscheme1}-\eqref{equ:admmscheme5}, we only give some key steps here.
	
	First, Lemma \ref{lem:descent} and \ref{lem:square} still work.  $\operatorname{dist}\left( 0,\partial\mathcal{L}_\beta(W^{k+1})\right) $ in Lemma \ref{lem:boundsub} is still bounded, but the upper bound is no longer as shown in \eqref{equ:boundsub}. Actually, \eqref{equ:Z} and \eqref{lem33-4} become
	\begin{equation*}
		\left\|Z^{k+1}-Z^{k+\frac12}\right\|_F\leq B_1\left\|X^{k+1}-X^k\right\|_F,
	\end{equation*}
	\begin{align*}
		\operatorname{dist}\left( 0,\partial\mathcal{L}_\beta(W^{k+1})\right) \le B_2&\Big( \left\|\Lambda^{k+1}-\Lambda^k\right\|_F+\left\|X^{k+1}-X^k\right\|_F \\
		&\hspace{0.5em} +\left\|Z^{k+1}-Z^k\right\|_F+\left\|Z^{k+1}-Z^{k+\frac12}\right\|_F\Big),
	\end{align*}
	respectively, where both $B_1$, $B_2>0$ are constants. Thus combining \eqref{equ:lambda}, \eqref{equ:X}, and the above two inequalities, there exists $\vartheta>0$ such that
	\begin{align*}
		\operatorname{dist}\left( 0,\partial\mathcal{L}_\beta(W^{k+1})\right) \le& \vartheta\left(\left\|Y^{k+1}-Y^k\right\|_F+\left\|Z^{k+1}-Z^k\right\|_F\right.\\
		&\hspace{1em}\left.+\left\|Y^{k}-Y^{k-1}\right\|_F+\left\|Z^{k}-Z^{k-1}\right\|_F\right).
	\end{align*}
	
	Lemma \ref{lem:subconverge} is consistent with the previous one. For the final convergence theorem \ref{thm:converge}, case 1 is unchanged. As for case 2, \eqref{thm1-4} becomes 
	\begin{align*}
		\frac{\Delta_{k,k+1}}{\varphi^{\prime}\left(\mathcal{L}_\beta(W^k)-\mathcal{L}_\beta(W^{\star})\right)} \le & \vartheta\left(\left\|Y^{k}-Y^{k-1}\right\|_F+\left\|Z^{k}-Z^{k-1}\right\|_F\right.\\
		&\hspace{1em}\left.+\left\|Y^{k-1}-Y^{k-2}\right\|_F+\left\|Z^{k-1}-Z^{k-2}\right\|_F\right)\cdot \Delta_{k,k+1}.
	\end{align*}
	Denote $A_k:=\left\|Y^{k+1}-Y^k\right\|_F+\left\|Z^{k+1}-Z^k\right\|_F$, then for any $k>k_1$,
	\begin{equation*}
		A_k^2\leq\frac{2\vartheta}c(A_{k-1}+A_{k-2})\cdot\Delta_{k,k+1}.
	\end{equation*}
	Thus,
	\begin{equation*}
		A_k\le \sqrt{\frac{4\vartheta}{c}\Delta_{k,k+1}}\cdot \sqrt{\frac12 \left( A_{k-1}+A_{k-2}\right) }\le \frac{2\vartheta}{c}\Delta_{k,k+1}+\frac14 \left( A_{k-1}+A_{k-2}\right) .
	\end{equation*}
	Adding $\frac{\sqrt{17}-1}{8}A_{k-1}$ to both sides of the inequality, we get
	\begin{equation*}
		A_k+\frac{\sqrt{17}-1}{8}A_{k-1}\le \frac{2\vartheta}{c}\Delta_{k,k+1}+\frac{\sqrt{17}+1}{8}\left(A_{k-1}+\frac{\sqrt{17}-1}{8}A_{k-2}\right).
	\end{equation*}
	Similar to \eqref{thm1-5}, we obtain 
	\begin{equation*}
		\sum_{k=0}^{\infty}A_k \le A_0+ \sum_{k=1}^{\infty}\left(A_k+\frac{\sqrt{17}-1}{8}A_{k-1}\right)<+\infty.
	\end{equation*}
	Proceeding as in the proof of Theorem \ref{thm:converge}, the final convergence is easy to be established.
\end{proof}

Another case is that practical problems do not satisfy Assumption \ref{ass:converge}: (1) $H(X,Y,Z)$ is nonseparable for $X$, $Y$, and $Z$, which means (a6) in Assumption \ref{ass:converge} is not satisfied. (2) The linear constraint in \eqref{equ:opt} does not contain either $X$ or $Y$, which implies one of $\mathcal{A}$ and $\mathcal{B}$ equals to the zero operator. Thus (a2) in Assumption \ref{ass:converge} is invalid. 

For case (1), we can apply the following proximal ADMMn to ensure the convergence:
\begin{align}\label{equ:2admmscheme}
	\left\{\begin{aligned}
		&X^{k+1}\in\mathop{\arg\min}\limits_{X} \left\{\mathcal{L}_{\beta}\left(X,Y^k,Z^k,\Lambda^k\right)+\frac12 \left\| X-X^k\right\|_Q^2\right\},  \\
		&Z^{k+\frac12}=\mathop{\arg\min}\limits_{Z} \mathcal{L}_{\beta}\left(X^{k+1},Y^k,Z,\Lambda^k\right),   \\
		&Y^{k+1}\in\mathop{\arg\min}\limits_{Y} \mathcal{L}_{\beta}\left(X^{k+1},Y,Z^{k+\frac12},\Lambda^k\right), \\
		&Z^{k+1}=\mathop{\arg\min}\limits_{Z} \mathcal{L}_{\beta}\left(X^{k+1},Y^{k+1},Z,\Lambda^k\right), \\
		&\Lambda^{k+1}=\Lambda^k-\beta\left(\mathcal{A}(X^{k+1})+\mathcal{B}(Y^{k+1})+\mathcal{C}(Z^{k+1})-D\right).
	\end{aligned}\right.
\end{align}
where $Q$ is a symmetric positive definite matrix, and $\|X\|_Q^2:=\left \langle QX,X \right \rangle$.

As for case (2), a possible approach is to rewrite an equivalent linear constraint so that it satisfies (a2) in Assumption \ref{ass:converge}. Another option is also adding a proximal term to the update of the missing variable in the linear constraint as in \eqref{equ:2admmscheme}. The convergence analysis of proximal ADMMn is easier than that of ADMMn \eqref{equ:admmscheme1}-\eqref{equ:admmscheme5} due to the existence of the proximal term, and can be deduced from the very recent work \cite{ConvDescentKL}. So we omit the proof here.

\section{Numerical experiments}\label{sec:4}

In this section, we make numerical experiments to test ADMMn for nonconvex optimization problems such as Robust Principal Component Analysis (RPCA) and Nonnegative Matrix Completion (NMC). For the sake of comparison, we also run PSR-ADMM \cite{PSRADMM}, pADMMz \cite{ztwice}, and IPPS-ADMM \cite{IPPSADMM} in the first experiment, and ADM \cite{ADM} for the Nonnegative Matrix Factorization/Completion (NMFC). All experiments are performed in MATLAB R2023a on a 64-bit laptop with an Intel Core i7-13700H 2.4GHz CPU and 16GB RAM.\footnote{All codes are available at \url{https://github.com/ZekunLiuOpt/ADMMn}.} 

\subsection{Robust principle component analysis}\label{sub:4-2}

Consider the following 3-block nonconvex model for RPCA \cite{PSRADMM,IPPSADMM}:
\begin{equation}\label{equ:rpca}
	\begin{aligned}
		\min_{X,Y,Z\in\mathbb{R}^{m\times n}}&\left\|X\right\|_*+ \rho \|Y\|_1 + \frac{\omega}{2}\|Z-M\|_F^2 \\ \mathrm{s.t.}\quad\quad& X+Y-Z=0,
	\end{aligned}
\end{equation}
where $M\in\R^{m\times n}$ is a given observation matrix, the low-rank term $\|X\|_*\!:=\! \sum_{i=1}^{\min \left\lbrace m,n\right\rbrace }\left|\sigma_i(X)\right|^{1/2}$ ($\sigma_i(X)$ is the singular value of $X$) and  $\|Y\|_1:=\sum_{i=1}^m\sum_{j=1}^n\left|Y_{ij}\right|$, $\rho$ and $\omega$ are two model parameters.

Note that \eqref{equ:rpca} is a form of problem \eqref{equ:opt} with $H\equiv0$. Applying the tested algorithms PSR-ADMM \cite{PSRADMM}, pADMMz \cite{ztwice}, IPPS-ADMM \cite{IPPSADMM}, and ADMMn to RPCA \eqref{equ:rpca}, the iteration details of ADMMn on RPCA  is 
\begin{equation*}
	\begin{cases}
		Y^{k+1}=\mathrm{prox}_{\frac{\rho}{\beta}\|\cdot\|_1}\left(Z^{k}+\frac{\Lambda^k}{\beta}-X^k\right), \\ 
		Z^{k+\frac12}=\frac{\omega M+\beta\left(X^k+Y^{k+1}\right)-\Lambda^k}{\omega+\beta}, \\
		X^{k+1}=\mathrm{prox}_{\frac{1}{\beta}\|\cdot\|_{1/2}^{1/2}}\left(Z^{k+\frac12}+\frac{\Lambda^k}{\beta}-Y^{k+1}\right),  \\
		Z^{k+1}=\frac{\omega M+\beta\left(X^{k+1}+Y^{k+1}\right)-\Lambda^k}{\omega+\beta}, \\
		\Lambda^{k+1}=\Lambda^k-\beta\left(X^{k+1}+Y^{k+1}-Z^{k+1}\right),
	\end{cases}
\end{equation*}
where $\mathrm{prox}_{\mu \|\cdot\|_1}$ is the soft shrinkage operator \cite{proxL1}, and $\mathrm{prox}_{\mu \|\cdot\|_{1/2}^{1/2}}$ is the half shrinkage operator \cite{proxL12}. Interested readers are referred to  \cite{PSRADMM,ztwice,IPPSADMM} for the update details of  PSR-ADMM, pADMMz and IPPS-ADMM on \eqref{equ:rpca}. 

Denote the sparsity rate of $Y$ and the rank of $X$ as ``spr'' and ``rank'', respectively. The numerical experiment setup in MATLAB is shown as follows:
\begin{lstlisting}
	X = randn(m,rank)*randn(rank,n); 
	q = randperm(m*n); Y = zeros(m,n);
	K = round(spr*m*n); Y(q(1:K)) = randn(K,1);
	sigma = 0; % Noiseless case; 
	sigma = 0.01; % Gaussian noise case; 
	N = randn(m,n)*sigma; Z = X + Y; M = Z + N;
\end{lstlisting}
Set $m=n=100$, $\rho=\frac{0.1}{\sqrt{m}}$, $\omega=10^3$, and test 8 different combinations of sparsity rate and rank with 20 independently random trials for each combination. $X$, $Y$, $Z$ and $\Lambda$ are all initialized to be zero matrices. For PSR-ADMM, set $F_1=F_2=0.07I$, $(r,s)=(-0.1,1.05)$ and $\beta_1=3.11$ as it suggested. Set $F_1=F_2=I$, $(r,s)=(0.2,0.5)$, $(\alpha_1,\alpha_2)=(0.25,0.25)$, $\beta_2=7.09$ for IPPS-ADMM as it recommended. As for ADMMn and pADMMz, we set $\beta_1=3.11$ to keep them consistent with PSR-ADMM, and set $S=T=\frac{I}{\beta_1}$ for pADMMz to keep its proximal terms consistent with IPPS-ADMM.
The stopping criterion is set to be
\begin{equation*}
	\mathrm{RelChg}:=\frac{\left\|(X^{k+1},Y^{k+1},Z^{k+1})-(X^k,Y^k,Z^k)\right\|_F}{\left\|(X^k,Y^k,Z^k)\right\|_F+1}\le\epsilon\mathrm{~or~}k>\mathrm{MaxIter},
\end{equation*}
where $\epsilon=10^{-7}$ is the tolerance and $\mathrm{MaxIter}=3000$.

Let the recovered solution of \eqref{equ:rpca} be $(\hat{X},\hat{Y},\hat{Z})$, and the ground truth is denoted as $(X^{\star},Y^{\star},Z^{\star})$. The relative error is used to measure the recovery quality:
\begin{equation*}
	\mathrm{RelErr}:=\frac{\left\|(\hat{X},\hat{Y},\hat{Z})-(X^{\star},Y^{\star},Z^{\star})\right\|_F}{\left\|(X^{\star},Y^{\star},Z^{\star})\right\|_F+1}.
\end{equation*}

\begin{figure}[htbp]
	\centering
	\begin{subfigure}{0.49\linewidth}
		\centering
		\includegraphics[width=1\linewidth]{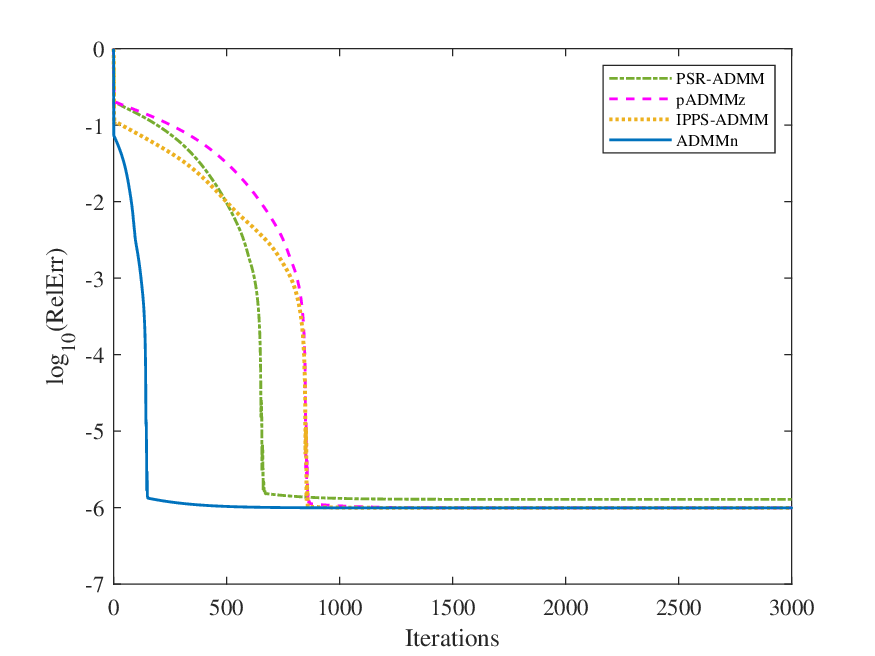}
		\caption{$\mathrm{spr} = 0.05, \mathrm{rank} = 10, \mathrm{noiseless}$}
		\label{f1}
	\end{subfigure}
	\centering
	\begin{subfigure}{0.49\linewidth}
	\centering
	\includegraphics[width=1\linewidth]{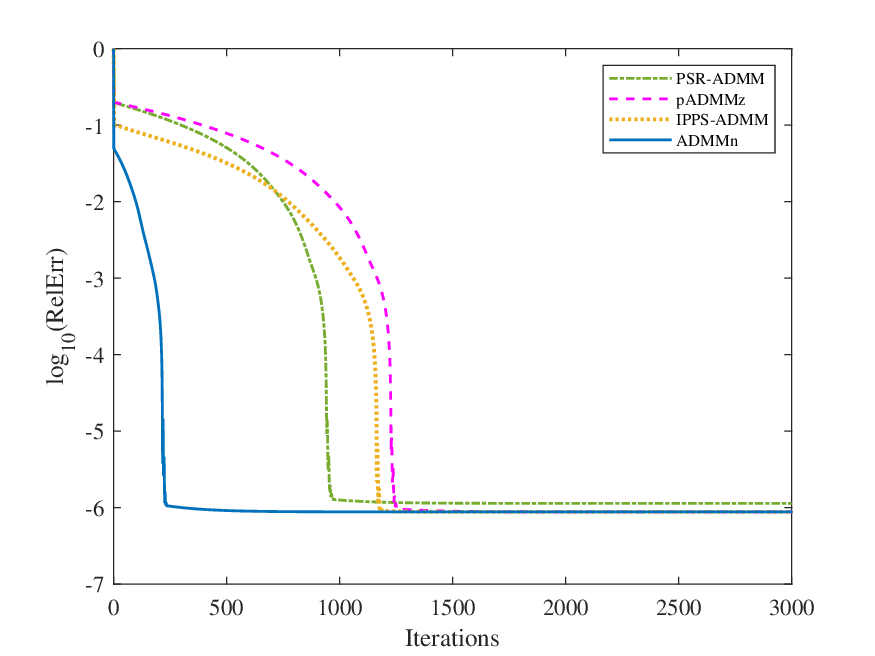}
	\caption{$\mathrm{spr} = 0.05, \mathrm{rank} = 20, \mathrm{noiseless}$}
	\label{f2}
	\end{subfigure}
	
	\centering
	\begin{subfigure}{0.49\linewidth}
		\centering
		\includegraphics[width=1\linewidth]{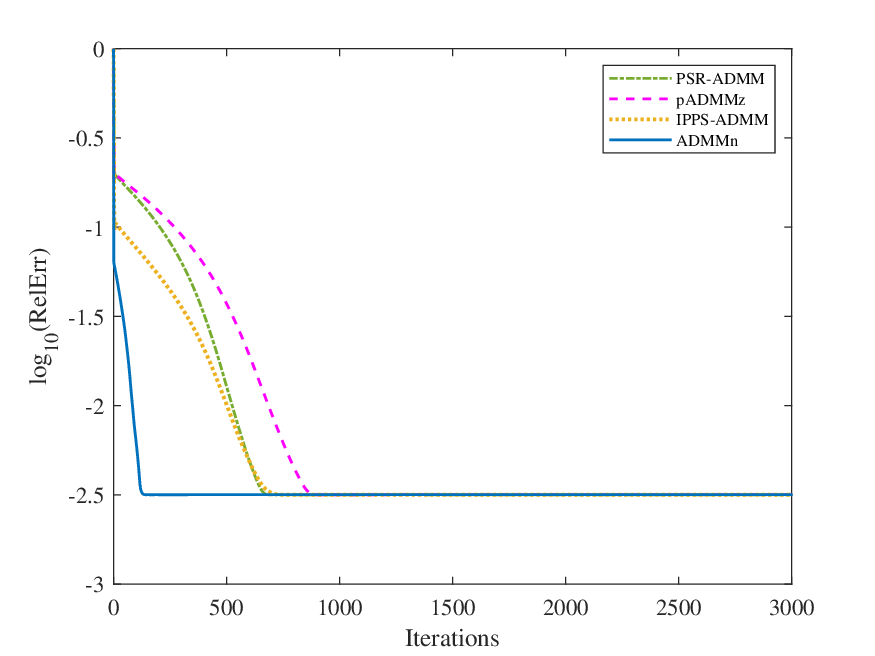}
		\caption{$\mathrm{spr} = 0.05, \mathrm{rank} = 10, \sigma=0.01$}
		\label{f3}
	\end{subfigure}
	\centering
	\begin{subfigure}{0.49\linewidth}
		\centering
		\includegraphics[width=1\linewidth]{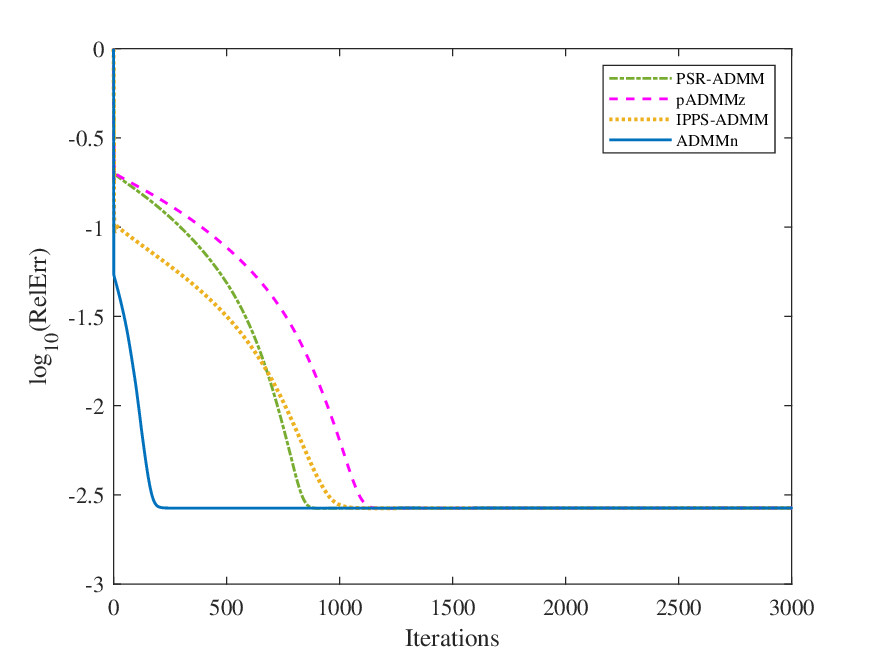}
		\caption{$\mathrm{spr} = 0.05, \mathrm{rank} = 20, \sigma=0.01$}
		\label{f4}
	\end{subfigure}
	\caption{The change of $\mathrm{log}_{10}(\mathrm{RelErr})$ with the iterations for the tested algorithms under different sparsity rate and rank.}
	\label{fig:rpca}
\end{figure}

Tables \ref{tab:rpcanoiseless} and \ref{tab:rpcanoisy} show the numerical results for the noiseless case and Gaussian noise $\sigma=0.01$ case, respectively, where iter represents the average number of iterations, and time is the average cpu time (seconds) to meet the stopping criterion. It can be seen that ADMMn outperforms the other tested algorithms on RPCA \eqref{equ:rpca}, since it can obtain a high precision solution with much less iterations and time for both noiseless and noisy cases. To see it much more clearly, we apply all the tested algorithms without stopping early (i.e., iterate to MaxIter) on (\ref{equ:rpca}). As shown in Figure \ref{fig:rpca}, ADMMn indeed needs much less iterations to converge.

\begin{table}[hbt!]
	\centering
	\resizebox{\textwidth}{!}{
		\begin{tabular}{cccccccccccccccccccccc}
			\hline 
			\multicolumn{3}{c}{(spr, rank)} & \multicolumn{5}{c}{PSR-ADMM} & \multicolumn{5}{c}{pADMMz} & \multicolumn{5}{c}{IPPS-ADMM}  & \multicolumn{4}{c}{ADMMn} \\
			\cline{1-2} \cline{4-7} \cline{9-12} \cline{14-17} \cline{19-22}
			spr & rank & & iter & time & RelChg & RelErr & & iter & time & RelChg & RelErr & & iter & time & RelChg & RelErr & & iter & time & RelChg & RelErr \\ \hline
			0.05 & 1 & & 362 & 0.32 & 5.3868e-08 & 4.2502e-06 & & 465 & 0.42 & 6.7297e-08 & 3.0163e-06 & & 387 & 0.34 & 6.2851e-08 & 2.7863e-06 & & 107 & 0.10 & 5.8242e-08 & 3.7968e-06 \\
			0.05 & 5 & & 529 & 0.47 & 7.7242e-08 & 1.9054e-06 & & 682 & 0.60 & 7.7340e-08 & 1.4082e-06 & & 567 & 0.50 & 5.7720e-08 & 1.3320e-06 & & 122 & 0.11 & 6.2023e-08 & 1.7386e-06 \\
			0.05 & 10 & & 718 & 0.62 & 7.4016e-08 & 1.4504e-06 & & 928 & 0.81 & 7.5172e-08 & 1.0847e-06 & & 783 & 0.68 & 6.9469e-08 & 1.0415e-06 & & 152 & 0.13 & 7.1687e-08 & 1.2837e-06 \\  
			0.05 & 20 & & 976 & 0.84 & 8.3950e-08 & 1.2983e-06 & & 1263 & 1.09 & 8.4113e-08 & 9.8996e-07 & & 1115 & 0.96 & 8.2786e-08 & 9.8948e-07 & & 214 & 0.18 & 7.8429e-08 & 1.1136e-06 \\
			0.1 & 1 & & 368 & 0.33 & 6.4783e-08 & 4.8912e-06 & & 473 & 0.42 & 6.2225e-08 & 3.5946e-06 & & 442 & 0.39 & 7.6569e-08 & 3.7810e-06 & & 163 & 0.15 & 6.2645e-08 & 4.4905e-06 \\
			0.1 & 5 & & 568 & 0.50 & 7.2070e-08 & 2.3899e-06 & & 732 & 0.65 & 7.0077e-08 & 1.8094e-06 & & 666 & 0.57 & 7.5995e-08 & 1.7467e-06 & & 193 & 0.17 & 6.8942e-08 & 2.0923e-06 \\
			0.1 & 10 & & 758 & 0.66 & 7.8605e-08 & 1.8989e-06 & & 978 & 0.85 & 8.3294e-08 & 1.4446e-06 & & 892 & 0.79 & 8.2104e-08 & 1.4166e-06 & & 238 & 0.21 & 7.0317e-08 & 1.5974e-06 \\
			0.1 & 20 & & 1113 & 0.97 & 8.9267e-08 & 1.7725e-06 & & 1435 & 1.25 & 8.7212e-08 & 1.3659e-06 & & 1461 & 1.28 & 8.2686e-08 & 1.3405e-06 & & 371 & 0.33 & 8.1308e-08 & 1.4093e-06 \\   
			\hline
		\end{tabular}
	}
	\caption{Comparison of the algorithms on RPCA with different sparsity rate and rank: $\sigma=0$.}
	\label{tab:rpcanoiseless}
\end{table}

\begin{table}[hbt!]
	\centering
	\resizebox{\textwidth}{!}{
		\begin{tabular}{cccccccccccccccccccccc}
			\hline 
			\multicolumn{3}{c}{(spr, rank)} & \multicolumn{5}{c}{PSR-ADMM} & \multicolumn{5}{c}{pADMMz} & \multicolumn{5}{c}{IPPS-ADMM}  & \multicolumn{4}{c}{ADMMn} \\
			\cline{1-2} \cline{4-7} \cline{9-12} \cline{14-17} \cline{19-22}
			spr & rank & & iter & time & RelChg & RelErr & & iter & time & RelChg & RelErr & & iter & time & RelChg & RelErr & & iter & time & RelChg & RelErr \\ \hline
			0.05 & 1 & & 1403 & 1.19 & 9.9078e-08 & 9.4924e-03 & & 1459 & 1.24 & 9.9404e-08 & 9.4945e-03 & & 1109 & 0.95 & 9.9427e-08 & 9.4945e-03 & & 1020 & 0.87 & 9.9185e-08 & 9.4944e-03 \\
			0.05 & 5 & & 1440 & 1.23 & 9.9784e-08 & 4.5857e-03 & & 1559 & 1.35 & 9.9790e-08 & 4.5866e-03 & & 1336 & 1.15 & 9.9803e-08 & 4.5866e-03 & & 934 & 0.81 & 9.9706e-08 & 4.5866e-03 \\
			0.05 & 10 & & 1583 & 1.33 & 9.9784e-08 & 3.3780e-03 & & 1777 & 1.50 & 9.9758e-08 & 3.3785e-03 & & 1618 & 1.36 & 9.9839e-08 & 3.3782e-03 & & 927 & 0.78 & 9.9746e-08 & 3.3786e-03 \\  
			0.05 & 20 & & 1868 & 1.58 & 9.9772e-08 & 2.6152e-03 & & 2175 & 1.84 & 9.9865e-08 & 2.6155e-03 & & 2128 & 1.80 & 9.9890e-08 & 2.6154e-03 & & 1033 & 0.87 & 9.9805e-08 & 2.6157e-03 \\
			0.1 & 1 & & 1405 & 1.26 & 9.9197e-08 & 9.4958e-03 & & 1456 & 1.30 & 9.9384e-08 & 9.4977e-03 & & 1197 & 1.07 & 9.9330e-08 & 9.4977e-03 & & 1078 & 0.96 & 9.9290e-08 & 9.4977e-03 \\
			0.1 & 5 & & 1496 & 1.33 & 9.9773e-08 & 4.6524e-03 & & 1633 & 1.45 & 9.9773e-08 & 4.6533e-03 & & 1484 & 1.31 & 9.9645e-08 & 4.6532e-03 & & 1018 & 0.91 & 9.9651e-08 & 4.6533e-03 \\
			0.1 & 10 & & 1650 & 1.43 & 9.9813e-08 & 3.4837e-03 & & 1871 & 1.63 & 9.9781e-08 & 3.4841e-03 & & 1803 & 1.56 & 9.9866e-08 & 3.4839e-03 & & 1045 & 0.91 & 9.9757e-08 & 3.4842e-03 \\
			0.1 & 20 & & 2084 & 1.78 & 9.9856e-08 & 2.8054e-03 & & 2455 & 2.10 & 9.9841e-08 & 2.8058e-03 & & 2610 & 2.23 & 1.0196e-07 & 2.8067e-03 & & 1254 & 1.07 & 9.9854e-08 & 2.8060e-03 \\   
			\hline
		\end{tabular}
	}
	\caption{Comparison of the algorithms on RPCA with different sparsity rate and rank: $\sigma=0.01$.}
	\label{tab:rpcanoisy}
\end{table}

\subsection{Nonnegative matrix completion}\label{sub:4-3}

The original rank-constrained nonnegative matrix completion problem can be formulated as
\begin{equation}\label{equ:originnmc}
	\begin{aligned}
		\min_{X\in\mathbb{R}^{m\times n}} &\left\|\mathcal{P}_{\Omega}(X-M)\right\|_F^2\\
		\mathrm{s.t.} \quad &\mathrm{rank}(X)\le r, \\
		\quad &X\ge 0,
	\end{aligned}
\end{equation}
where $M$ is the observation matrix, $r$ is the given upper rank estimation of the matrix $X$, and $\mathcal{P}_{\Omega}$ is the projection onto the sampling set $\Omega$:
\begin{align*}
	\left(\mathcal{P}_{\Omega}(X)\right)_{ij}=
	\begin{cases}~X_{ij},\quad&\text{if}~(i,j)\in\Omega,\\
		~0,\quad&\text{if}~(i,j)\notin\Omega.
	\end{cases}
\end{align*}

If we reformulate the above model to a 2-block nonconvex problem with a linear constraint, either the subproblem is hard to solve, or none of existed methods on the reformulated problem have convergence guarantee. Hence we model the above problem to a 3-block nonconvex form as follows:
\begin{equation}\label{equ:nmc}
	\begin{aligned}
		\min_{X,Y,Z\in\mathbb{R}^{m\times n}} &\delta_{\mathcal{K}}(X)+\delta_{\mathcal{N}}(Y)+\left\|\mathcal{P}_{\Omega}(Z-M)\right\|_F^2+\frac\rho2 \|Y-Z\|_F^2\\
		\mathrm{s.t.} \quad\quad &2X-Y-Z=0,
	\end{aligned}
\end{equation}
where $\delta_{\mathcal{K}}$, $\delta_{\mathcal{N}}$ are indicator functions to $\mathcal{K}:=\left\{X\in\mathbb{R}^{m\times n}:\mathrm{rank}(X)\le r\right\}$, $\mathcal{N}:=\left\{X\in\mathbb{R}^{m\times n}:X\ge 0\right\}$, respectively, and $\rho>0$ is the penalty parameter.

Note that \eqref{equ:nmc} is a form of \eqref{equ:opt} with the nonseparable term $H:=\frac\rho2 \|Y-Z\|_F^2$. Since the algorithms in the first experiment are not designed for nonseparable problems like \eqref{equ:nmc}, we only apply ADMMn to NMC \eqref{equ:nmc}. The update format of ADMMn is
\begin{equation*}
	\begin{cases}
		Y^{k+1}=\frac{1}{\rho+\beta}\max\left(2\beta X^{k}+(\rho-\beta) Z^{k}-\Lambda^k,0\right), \\ 
		Z^{k+\frac12}=\left(2\Omega \odot M+(\rho-\beta) Y^{k+1}+2\beta X^k-\Lambda^k\right)\oslash \left(2\Omega+(\rho+\beta) \textbf{1}\right), \\
		X^{k+1}=\frac12 \mathrm{TSVD}\left(Y^{k+1}+Z^{k+\frac12}+\frac{\Lambda^k}{\beta},r\right),  \\
		Z^{k+1}=\left(2\Omega \odot M+(\rho-\beta) Y^{k+1}+2\beta X^{k+1}-\Lambda^k\right)\oslash \left(2\Omega+(\rho+\beta) \textbf{1}\right), \\
		\Lambda^{k+1}=\Lambda^k-\beta\left(2X^{k+1}-Y^{k+1}-Z^{k+1}\right),
	\end{cases}
\end{equation*}
where $\mathrm{TSVD}(\cdot,r)$ is the Truncated Singular Value Decomposition (TSVD) \cite{TSVD}, and $\textbf{1}$ represents the $m\times n$ matrix whose entries all equal to 1. 

For comparison, we consider another formulation of \eqref{equ:originnmc} which is known as the Nonnegative Factorization Matrix Completion (NFMC) problem as follows:
\begin{equation}\label{equ:nfmc}
	\begin{aligned}
		\operatorname*{min}_{U,V,X,Y,Z}\quad & \frac{1}{2}\|XY-Z\|_{F}^{2}   \\ 
		\mathrm{s.t.} \quad\quad&X=U,\; Y=V,  \\
		&U\ge0,\; V\ge0, \\
		&\mathcal{P}_{\Omega}(Z-M)=0,
	\end{aligned}
\end{equation}
where $X,U\in\mathbb{R}^{m\times r}$ and $Y,V\in\mathbb{R}^{r\times n}$. An ADMM scheme algorithm is proposed in \cite{ADM} to solve the 5-block problem \eqref{equ:nfmc} and has the global convergence guarantee. We call it ADM in this experiment and compare ADMMn with it. 

Let ``rank'' and ``sr'' represent the rank $r$ of the original low-rank matrix and the sampling rate $|\Omega|/mn$, respectively. The MATLAB code for experiment setup is shown as follows:
\begin{lstlisting}
	U = rand(m,rank); V = rand(rank,n); X = U*V; 
	% Generate the low-rank matrix 
	Omega = zeros(m,n); 
	ind = randi(numel(Omega),1,round(m*n*sr));
	Omega(ind) = 1; % Generate the sampling set
	M = Omega.*X; % The observation matrix   
\end{lstlisting}
Set $m=n=500$, and test 9 combinations of rank and sampling rate with 20 independently random trials for each combination. For ADMMn on \eqref{equ:nmc}, set $\rho=1$, $\beta_1=1$, and $X$, $Y$, $Z$, $\Lambda$ are all initialized as zero matrices. To keep consistent with the settings of ADM on \eqref{equ:nfmc}, we set the stopping criterion as
\begin{equation*}
	\mathrm{RelChg}:=\frac{\left\|\mathcal{P}_{\Omega}(X^k-M)\right\|_F}{\|M\|_F}\le\epsilon\mathrm{~or~}k>\mathrm{MaxIter},
\end{equation*}
where $\epsilon=10^{-6}$ is the tolerance and $\mathrm{MaxIter}=3000$.

Denote the recovered solution of \eqref{equ:nmc} as $(\hat{X},\hat{Y},\hat{Z})$, and the ground truth is $(X^{\star},Y^{\star},Z^{\star})$. We use the following relative error to measure the recovery quality in order to keep consistent with the ADM on \eqref{equ:nfmc}:
\begin{equation*}
	\mathrm{RelErr}:=\frac{\left\|\hat{X}-X^{\star}\right\|_F}{\left\|X^{\star}\right\|_F+1}.
\end{equation*}

For ADM on \eqref{equ:nfmc}, as \cite{ADM} suggested, set $\gamma=1.618$, $\alpha=1.91\times 10^{-4}\|M\|_F\frac{\max \left\lbrace m,n\right\rbrace }{r}$ and $\beta_2=\frac{n\alpha}{m}$. $U$, $V$, $\Lambda$, $\Pi$ are initialized as zero matrices, while $Y$ is set as a nonnegative random matrix, and $Z=M$. The stopping criterion in \cite{ADM} is given as
\begin{equation*}
	\mathrm{RelChg}:=\frac{\left\|\mathcal{P}_{\Omega}(X^kY^k-M)\right\|_F}{\left\|M\right\|_F}\le\epsilon\mathrm{~or~}k>\mathrm{MaxIter},
\end{equation*}
where $\epsilon=10^{-6}$ is the tolerance and $\mathrm{MaxIter}=3000$.

Denote the recovered solution of \eqref{equ:nfmc} as $(\hat{X},\hat{Y},\hat{Z},\hat{U},\hat{V})$, and the ground truth is $(X^{\star},Y^{\star},Z^{\star},U^{\star},V^{\star})$. \cite{ADM} employs the relative error to measure the recovery quality:
\begin{equation*}
	\mathrm{RelErr}:=\frac{\left\|\hat{X}\hat{Y}-X^{\star}Y^{\star}\right\|_F}{\left\|X^{\star}Y^{\star}\right\|_F+1}.
\end{equation*}

Table \ref{tab:nmc} shows the experiment results. We can see that when the original matrix is really low-rank, both ADMMn and ADM can recover it with high precision, and ADMMn uses much less iterations and time in reconstruction than ADM. However, when the rank of the original matrix is a little large, the recovery quality of ADM is far inferior to ADMMn, and the stopping criterion of ADM is not satisfied in the iterations. Besides, we can find that ADMMn would cost more time than ADM when their numbers of iterations are the same. This occurs because TSVD is applied in NMC \eqref{equ:nmc}, which is more time-consuming than the matrix factorization in NFMC \eqref{equ:nfmc}. Benefited from the fewer iterations required to converge, ADMMn spends less time than ADM in general.  

\begin{table}[hbt!]
	\centering
	\resizebox{\textwidth}{!}{
		\begin{tabular}{cccccccccccc}
			\hline 
			\multicolumn{3}{c}{(rank, sr)} & \multicolumn{5}{c}{ADM} & \multicolumn{4}{c}{ADMMn} \\
			\cline{1-2} \cline{4-7} \cline{9-12}
			rank & sr & & iter & time & RelChg & RelErr & & iter & time & RelChg & RelErr \\ \hline
			2 & 0.7 & & 2833 & 3.52 & 9.9568e-06 & 9.9090e-06 & & 93 & 0.30 & 9.4517e-07 & 1.0523e-06 \\ 
			2 & 0.5 & & 2806 & 3.48 & 1.5677e-05 & 1.5935e-05 & & 125 & 0.38 & 9.6058e-07 & 1.1160e-06 \\ 
			2 & 0.3 & & 2966 & 3.71 & 3.6358e-05 & 3.7457e-05 & & 214 & 0.62 & 9.6896e-07 & 1.2213e-06 \\ 
			10 & 0.7 & & 3000 & 4.13 & 3.5762e-04 & 3.7393e-04 & & 130 & 0.46 & 9.4706e-07 & 1.2106e-06 \\ 
			10 & 0.5 & & 3000 & 4.07 & 7.8137e-04 & 8.4139e-04 & & 197 & 0.67 & 9.6856e-07 & 1.3378e-06 \\ 
			10 & 0.3 & & 3000 & 4.07 & 4.5016e-04 & 5.2008e-04 & & 414 & 1.41 & 9.8436e-07 & 1.5899e-06 \\ 
			20 & 0.7 & & 3000 & 4.16 & 2.6131e-03 & 2.8698e-03 & & 223 & 1.19 & 9.6093e-07 & 1.3835e-06 \\ 
			20 & 0.5 & & 3000 & 4.16 & 2.9042e-03 & 3.3680e-03 & & 376 & 2.00 & 9.7492e-07 & 1.5764e-06 \\ 
			20 & 0.3 & & 3000 & 4.17 & 3.3527e-03 & 4.3619e-03 & & 1186 & 6.29 & 9.8860e-07 & 2.0521e-06 \\ 
			\hline
		\end{tabular}
	}
	\caption{Comparison of the algorithms on NMC/NFMC with different rank and sampling rate.}
	\label{tab:nmc}
\end{table}

\section{Conclusions}\label{sec:5}

In summary, we propose an ADMM algorithm with the third variable updated twice at each iteration to solve the 3-block nonconvex nonseparable problem with a linear constraint. Global convergence of the proposed ADMM is established under the Kurdyka-Łojasiewicz property. Besides, we discuss two simple extensions of ADMMn which are useful in applications. At last, we experiment on a 3-block separable nonconvex RPCA problem with $H\equiv0$, and a brand new 3-block nonseparable nonconvex NMC problem proposed by us, respectively. The numerical results show that for both \eqref{equ:opt} and its degeneration form, the proposed ADMMn is consistent with theoretical results and exhibits a clear advantage in practice. \\

%\textbf{Acknowledgments.} The author would like to thank his supervisor, Prof. Jinyan Fan, for her insightful and constructive comments that significantly improves the paper.

\noindent \textbf{Code Availability.} The source code is made available and can be obtained from \url{https://github.com/ZekunLiuOpt/ADMMn}.

	\end{document}